\documentclass[a4paper]{article}

\usepackage{amsmath}
\usepackage{framed}
\usepackage[utf8x]{inputenc}
\usepackage[dvipsnames]{xcolor}
\usepackage{xspace}
\usepackage{graphicx}
\usepackage{amsfonts} 
\usepackage{amssymb}
\usepackage{xcolor,cancel}
\usepackage{url}
\usepackage{rotating}
\usepackage{authblk}
\usepackage{subfigure}

\DeclareMathOperator{\argmin}{argmin}
\DeclareMathOperator{\supp}{supp}
\DeclareMathOperator{\Span}{Span}
\DeclareMathOperator{\sign}{sign}

\newcommand{\uargmin}[1]{\underset{#1}{\argmin}\;}
\newcommand{\norm}[1]{\| #1 \|}
\newcommand{\ie}{{\em i.e.,~}}
\newcommand{\eg}{{\em e.g.,~}}
\newcommand{\lcf}{{\em cf.~}}
\newcommand{\prox}{\mathrm{prox}}
\newcommand{\enscond}[2]{ \left\{ #1 \;:\; #2 \right\} }
\newcommand{\Ii}{\mathcal{I}}
\newcommand{\ins}[1]{\mathrm{#1}}
\newcommand{\dotp}[2]{\left\langle #1,\,#2 \right\rangle}
\newcommand{\Id}{\ins{Id}}
\newcommand{\qandq}{ \quad \text{and} \quad }
\newcommand{\zt}{\tilde z}
\newcommand{\xt}{\tilde x}
\newcommand{\x}{x}
\newcommand{\choice}[1]{ %
	\left\{ %
		\begin{array}{ll} #1 \end{array} %
	\right. }
\newcommand{\ifq}{{\text{if} \quad}}

\newcommand{\otherwise}{{\text{otherwise}}}
\newcommand{\qwhereq}{ \quad \text{where} \quad }
\newcommand{\andq}{ \text{and} \quad }

\newtheorem{defn}{Definition}

\usepackage{geometry} 
\newcommand{\RR}{\mathbb{R}}

\newcommand{\snorm}[1]{\left\|#1\right\|}

\begin{document}
\title{Refitting  solutions promoted by $\ell_{12}$ sparse analysis regularization with block penalties}

\author[1,2]{Charles-Alban Deledalle}
\author[2]{Nicolas Papadakis}
\author[3]{Joseph Salmon}
\author[4]{Samuel Vaiter}
\date{}      
%
\affil[1]{University  of  California,   San  Diego,   La  Jolla,   USA}
\affil[2]{CNRS, Univ. Bordeaux, IMB, F-33400 Talence, France}

\affil[3]{IMAG, Univ.   Montpellier,    CNRS,   Montpellier,    France}
 
\affil[4]{IMB, CNRS, Universit\'e de Bourgogne, 21078 Dijon, France}
\maketitle    

\begin{abstract}
In inverse problems, the use of an $\ell_{12}$ analysis regularizer induces a bias in the estimated solution.
We propose a general refitting framework for removing this artifact while keeping information of interest contained in the biased solution. This is done through the use of refitting block penalties that only act on the co-support of the estimation. Based on an analysis of related works in the literature, we propose a new penalty that is well suited for refitting purposes.
We also present an efficient algorithmic method to obtain the refitted solution along with the original (biased) solution for any convex refitting block penalty. Experiments illustrate the good behavior of the proposed block penalty for refitting.

\end{abstract}

\section{Introduction}

We consider linear inverse problems of the form $y = \Phi x + w$,
where $y \in \RR^p$ is an observed degraded image, $x \in \RR^n$ the unknown clean image, $\Phi: \RR^{n} \to \RR^p$
a linear operator and $w \in \RR^p$ a noise component, typically a zero-mean
white Gaussian random vector with standard deviation $\sigma > 0$.
To reduce the effect of noise and the potential ill-conditioning of $\Phi$,
we consider a regularized least square problem with an $\ell_{12}$ structured sparse analysis term of the form
\begin{align}\label{isotv}
  \hat{x} \in
  \uargmin{x} \tfrac12 \norm{\Phi x - y }^2 + \lambda \norm{\Gamma x}_{1,2}~.
\end{align}
where $\lambda > 0$ is a regularization parameter, $\Gamma : \RR^{n} \to \RR^{m \times b}$ is a linear analysis operator
mapping an image over $m$ blocks of size $b$ and $\norm{ z}_{1,2}=\sum_{i=1}^m \norm{z_i}=\sum_{i=1}^m( \sum_{j=1}^b z_{i,j}^2)^{1/2}$, with $z_i=\{z_{i,j}\}_{j=1}^b\in\RR ^b$.
This model is known to recover co-sparse solutions, \ie such that
$(\Gamma x)_i = 0_b$ for most blocks $1 \leq i \leq m$.
A typical example is the one of isotropic total-variation (TViso)
with $\Gamma = \nabla$ being the operator which extracts $m=n$
image gradient vectors of size $b = 2$ (for volumes $b = 3$, and so on).
The anisotropic total-variation is another example corresponding to
$\Gamma$ the operator which concatenates the vertical and horizontal
components of the gradients into a vector of size $m=2n$, hence $b=1$.

\subsection{Refitting}

\begin{figure}[!ht]
  \centering
  \vspace{-1em}
  \subfigure{\includegraphics[width=0.137\linewidth]{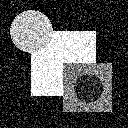}}\hfill%
  \subfigure{\includegraphics[width=0.137\linewidth]{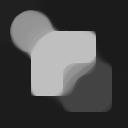}}\hfill%
    \subfigure{\includegraphics[width=0.137\linewidth]{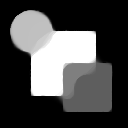}}\hfill%
  \subfigure{\includegraphics[width=0.137\linewidth]{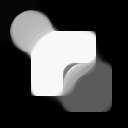}}\hfill%
  \subfigure{\includegraphics[width=0.137\linewidth]{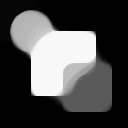}}\hfill%
  \subfigure{\includegraphics[width=0.137\linewidth]{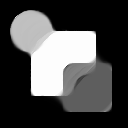}}\hfill%
  \subfigure{\includegraphics[width=0.137\linewidth]{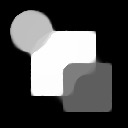}}\\[-1em]%
  \setcounter{subfigure}{0}%
  \subfigure[Noisy $y$]{\includegraphics[width=0.137\linewidth,viewport=94 84 222 212, clip]{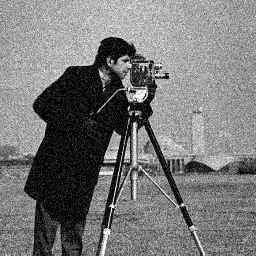}}\hfill%
  \subfigure[TViso $\hat x$]{\includegraphics[width=0.137\linewidth,viewport=94 84 222 212, clip]{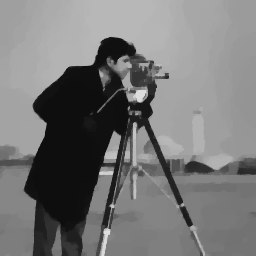}}\hfill%
    \subfigure[IB \cite{Osher}]{\includegraphics[width=0.137\linewidth,viewport=94 84 222 212, clip]{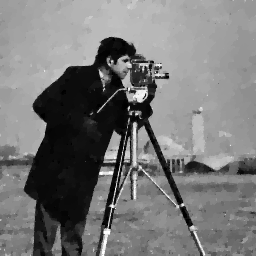}}\hfill%
  \subfigure[HO \cite{brinkmann2016bias}]{\includegraphics[width=0.137\linewidth,viewport=94 84 222 212, clip]{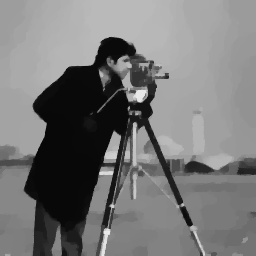}}\hfill%
  \subfigure[HD \cite{brinkmann2016bias}]{\includegraphics[width=0.137\linewidth,viewport=94 84 222 212, clip]{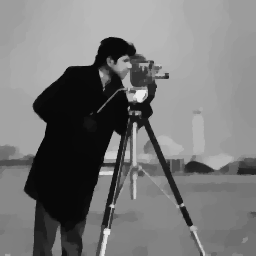}}\hfill%
  \subfigure[QO \cite{deledalle2016clear}]{\includegraphics[width=0.137\linewidth,viewport=94 84 222 212, clip]{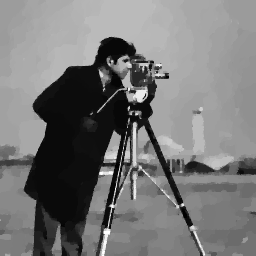}}\hfill%
  \subfigure[New SD]{\includegraphics[width=0.137\linewidth,viewport=94 84 222 212, clip]{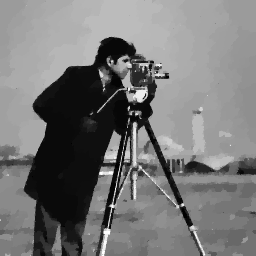}}\vspace*{-.9em}\\
  \caption{\label{ex:TViso} Comparison of standard refitting approaches with the proposed SD model.}
  \vspace*{-0em}
\end{figure}

The co-support of an image $x$ (or support of $\Gamma x$) is the set of its non-zero blocks:
\begin{align}
  \supp(\Gamma x) = \enscond{1 \leq i \leq m}{(\Gamma x)_i \ne 0_b}~.
\end{align}
While in some cases, the estimate $\hat{x}$ obtained by structured sparse analysis regularization \eqref{isotv}
recovers correctly the co-support of the underlying signal $x$,
it nevertheless suffers from a systematical bias in the estimated amplitudes $\hat{x}_i$.
With TViso, this bias is reflected by a loss of contrast
(see Fig.~\ref{ex:TViso}(b)).
A standard strategy to reduce this effect, called refitting \cite{efron2004least,Rigollet_Tsybakov11,Belloni_Chernozhukov13,Lederer13}, consists in approximating $y$ through $\Phi$ 
by an image sharing the same co-support as $\hat{x}$:
\begin{equation}\label{proj_supp}
  \tilde{x}^{\supp} \in \uargmin{x; \; \supp(\Gamma x) \subseteq \hat{\Ii}}
  \tfrac12 \| \Phi x - y \|^2 \enspace,\vspace{-0.1cm}
\end{equation}
where $\hat{\Ii} = \supp(\Gamma \hat{x})$.
While this strategy works well for blocks of size $b=1$, it suffers from an excessive increase of variance whenever $b\geq 2$, \eg for TViso.
This is due to the fact that solutions do not only present sharp edges,
but may involve gradual transitions.
To cope with this issue, additional features of $\hat x$ than its co-support must be preserved by a refitting procedure.
For the LASSO ($\Gamma=\Id$, $m=n$ and $b=1$), a pointwise preservation of the sign of $\hat x_i$ onto the support improves the numerical performances of the refitting \cite{2017arXiv170705232C}.

\subsection{Outline and contributions} In this paper, we introduce a new framework for refitting solutions promoted by $\ell_{12}$ structured sparse analysis regularization \eqref{isotv}. In Section 2, we present related works, illustrated in Figure \ref{ex:TViso}, that include Bregman iterations\cite{Osher} or debiasing approaches \cite{deledalle2016clear,brinkmann2016bias}. 
In Section 3, we describe our general variational refitting method for block penalties and show how the works presented in  Section 2 can be described with such a framework.
We  discuss suitable properties a refitting block penalty should satisfy and introduce the Soft-penalized Direction model (SD), a new flexible refitting block penalty inheriting the advantages of the Bregman-based approaches.
In Section 4, we propose a stable and one-step algorithm to compute our refitting strategy for any convex refitting block penalty, including the models in \cite{brinkmann2016bias} and \cite{deledalle2016clear}.
Experiments in Section 5 illustrate the practical benefits for the SD refitting and the potential of our framework for image processing.

\section{Related re-fitting works}
We first present some properties of Bregman divergences  used all along the paper.
\subsection{Properties of Bregman divergence of $\ell_{12}$ structured regularizers}\label{sec:intro_breg}
In the literature \cite{Osher,brinkmann2016bias}, Bregman divergences have proven to be well suited to measure the discrepancy between the biased solution $\hat x$ and its refitting $\tilde x$.
We recall that for a convex function $\psi$, the associated (generalized) Bregman divergence between $x$ and $\hat x$ is, for any subgradient $\hat p\in \partial \psi(\hat x)$: 
$D_\psi(x,\hat x)=\psi(x)-\psi(\hat x)-\langle \hat p,x-\hat x\rangle\geq 0\enspace.$
If $\psi$ is an absolutely 1-homogeneous function ($\psi(\alpha x)=|\alpha| \psi(x)$, $\forall\alpha \in\RR$) then $p\in \partial \psi(x)\Rightarrow \psi(x)=\dotp{p}{x}$ and the Bregman divergence simplifies into $D_\psi(x,\hat x)=\psi(x)-\langle \hat p,x\rangle\geq 0$.
For regularizers of the form $\psi(x)=\|\Gamma x\|_{1,2}$, a subgradient $\hat p\in \partial\|\Gamma \cdot\|_{1,2}(\hat x)$ satisfies $ \hat p=\Gamma^\top  \hat z$  with \cite{BGMEC16}: 
\begin{equation}\label{opt_z}
\left\{\begin{array}{ll}
\hat z_i=\tfrac{(\Gamma \hat x)_i}{\|(\Gamma \hat x)_i\|}\quad \quad & \textrm{if } i\in \hat{\Ii} =\supp(\Gamma \hat x)\enspace,\\
\norm{\hat z_i}\leq 1&\textrm{otherwise} \enspace,
\end{array}\right.
\end{equation}%
where we have
$\Gamma^\top\hat  z\in\partial \norm{\Gamma \cdot}_{1,2}(\hat x)\Leftrightarrow \hat z_i\in\partial \norm{\cdot}((\Gamma \hat x)_i),\, \forall i \in [m].$
Now denoting 
\begin{equation}\label{D_i}
D_i(\Gamma x)=\norm{(\Gamma x)_i} -\langle  \hat  z_i, (\Gamma x)_i\rangle=D_{ \norm{\cdot}} ((\Gamma x)_i,(\Gamma \hat x)_i),
\end{equation}
there exists an interesting relation between the global Bregman divergence on  $x$'s  and the  local ones on $(\Gamma x)_i$'s:\vspace{-0.1cm}
\begin{equation}\label{breg}
D_{\norm{\Gamma \cdot}_{1,2}} (x,\hat x)
=\sum_{i=1}^m \left( \norm{(\Gamma x)_i} -\langle  \hat  z_i, (\Gamma x)_i\rangle\right)=
\sum_{i=1}^m D_i(\Gamma x)
\enspace.\vspace{-0.2cm}
\end{equation}
Combining relations \eqref{opt_z} and \eqref{D_i}, 
we see that on the co-support $i\in\hat{\Ii}$,  such divergence  measures the fit of directions between $(\Gamma x)_i$ and $(\Gamma \hat x)_i$:
\begin{equation}\forall i\in \hat{\Ii},\, D_i(\Gamma x)=0\textrm{ iff } \exists \alpha_i\geq 0\textrm{ such that }(\Gamma  x)_i=\alpha_i (\Gamma \hat x)_i \enspace.\label{sub_supp1}\end{equation}
This divergence also partially captures the co-support as we have from  \eqref{D_i}
\begin{equation}
  \forall i\in \hat{\Ii}^c, \,
  D_i(\Gamma x) = 0 
\textrm{ with }
 \hat z_i\in \partial \norm{\cdot}((\Gamma \hat x)_i)\textrm{ s.t. }\hspace{-0.03cm}\norm{\hat z_i}<1
\textrm{ iff }  {(\Gamma x)_i}=0_b\enspace.\label{sub_supp2}
\end{equation}

\subsection{Bregman-based Refitting}
We now review some refitting methods based on Bregman divergences.
\vspace{-.5em}
\subsubsection{Flexible Iterative Bregman regularization}
The  Bregman process \cite{Osher} reduces the bias of solutions of \eqref{isotv}
by successively solving problems of the form
\begin{align}\label{iterative_OSHER}
  &
  \tilde{x}_{l+1} \in\uargmin{x}
  \tfrac12 \|\Phi x - y \|^2 + \lambda  D_{\norm{\Gamma \cdot}_{1,2}} (x,\tilde{x}_{l}).
\end{align}
We here consider a fixed $\lambda$, but different refitting strategies can be considered with decreasing parameters $\lambda_l$ as in \cite{Schetzer01,Tadmor04}.

\noindent
Setting  $\tilde{x}_{0}=0_n$, we have  $0_n\in \partial \norm{\Gamma \cdot}_{1,2}(0_n)$ so that
$D_{\norm{\Gamma \cdot}_{1,2}} (x,0_n)= \norm{\Gamma x}_{1,2}$ and we recover the biased solution of \eqref{isotv} as $\tilde{x}_{1}=\hat x$. 
We denote by $\tilde{x}^{\mathrm{IB}}=\tilde{x}_{2}$ the refitting obtained after $2$ steps of the Iterative Bregman (IB)  procedure  \eqref{iterative_OSHER}:
\begin{align}\label{iterative_refit}\tilde{x}^{\mathrm{IB}}=\tilde{x}_{2}\in\uargmin{x} \tfrac12 \|\Phi x - y \|^2 +\lambda D_{\norm{\Gamma \cdot}_{1,2}} (x,\hat x)\enspace.
\end{align}
As underlined in relation \eqref{sub_supp1}, by minimizing $D_{\norm{\Gamma \cdot}_{1,2}} (x,\hat x)$, one aims at preserving the direction of $\Gamma \hat x$ on the support $\hat{\Ii}$, without ensuring $\supp(\Gamma \tilde x^\textrm{IB})\subseteq \hat{\Ii}$.
This issue can be observed in the background of {\it Cameraman} in  Fig.~\ref{ex:TViso}(c), where noise reinjection is visible.
For the iterative framework, the co-support of the previous solution may indeed not be preserved ($\|(\Gamma \tilde x_l)_i\|=0 \nRightarrow \|(\Gamma \tilde x_{l+1})_i\|=0$) and can hence grow. The support of $\Gamma x_0$ for $\tilde{x}_{0}=0_n$ is for instance totally empty whereas the one of  $\hat x= \tilde{x}_{1}$ may not  (and should not) be empty.
For $l\to\infty$, the process actually converges to some $x$ such that $\Phi x=y$.

\vspace{-.5em}
\subsubsection{Hard-constrained refitting without explicit support identification.}
In order to respect the support of the biased solution $\hat x$ and to keep track of the direction $\Gamma \hat x$ during the refitting, the authors of \cite{brinkmann2016bias} proposed the following model:
\begin{align}\label{refit_german:HD}
  &
  \tilde{x}^{\mathrm{HD}} \in\uargmin{x; \hat  p\in \partial \norm{\Gamma \cdot}_{1,2}(x)}
  \tfrac12 \| \Phi x - y \|^2 \enspace,
\end{align}
for  $\hat p\in \partial \norm{\Gamma \cdot}_{1,2}(\hat x)$.
This model enforces  the  Bregman divergence  to be  $0$,
since  $\hat  p\in \partial \norm{\Gamma \cdot}_{1,2}(x)\Rightarrow \norm{\Gamma x}_{1,2}=\langle \hat p, x\rangle \Rightarrow D_{\norm{\Gamma \cdot}_{1,2}}(x,\hat x)=0$.

We see from \eqref{sub_supp1} that for $i\in\hat{\Ii}$, the direction of $(\Gamma \hat x)_i$ is preserved in the refitted solution. Following relations \eqref{opt_z} and \eqref{sub_supp2}, the co-support is also preserved for any $i\in \hat{\Ii}^c$ such that $\norm{\hat z_i}<1$.
Note though that extra elements in the co-support $\Gamma \tilde x^\textrm{HD}$ may be added at coordinates $i\in\hat{\Ii}^c$ such that $\norm{\hat{z}_i}=1$.
We denote this model as HD, for Hard-constrained Direction. To get ride of the direction dependency, a Hard-constrained Orientation (HO) model is also proposed in \cite{brinkmann2016bias}:
\begin{align}\label{refit_german:HO}
  &
  \tilde{x}^{\mathrm{HO}} \in\uargmin{x; \pm \hat  p\in \partial \norm{\Gamma \cdot}_{1,2}(x)}
  \tfrac12 \| \Phi x - y \|^2 \enspace.
\end{align}
The orientation model may nevertheless involve contrast inversions between biased and refitted solutions, as shown in Fig.~\ref{ex:TViso}(d) with the banana dark shape in the white region.
In practice, relaxations 
are used in \cite{brinkmann2016bias}
by solving, for a large value $\gamma>0$
\begin{align}\label{refit_german:HD2}
  &
  \tilde{x}_\gamma^{\mathrm{HD}} \in\uargmin{x}
  \tfrac12 \| \Phi x - y \|^2 +\gamma D_{\norm{\Gamma \cdot}_{1,2}} (x,\hat x)\enspace.
\end{align}
The main advantage of this refitting strategy is that no support identification is required since everything is implicitly encoded in the subgradient $\hat p$. This makes the process stable even if  the estimation of $\hat x$ is not highly accurate.
The support of $\Gamma \hat x$ is nevertheless only approximately preserved, since the constraint $D_{\norm{\Gamma \cdot}_{1,2}} (x,\hat x)=0$ can never be ensured numerically with a finite value of $\gamma$.
Finally, as shown in Fig.~\ref{ex:TViso}(d-e), such constrained approaches  lack of flexibility since the  orientation of $\Gamma \tilde x$ cannot deviate from the one of $\Gamma \hat x$
(for complex signals, such as {\it Cameraman}, amplitudes remain significantly biased and less details are recovered).

\subsection{Flexible Quadratic refitting without support identification}

We now describe an alternative way for performing variational refitting.
When specialized to $\ell_{1,2}$ sparse analysis regularization, CLEAR, a general refitting framework \cite{deledalle2016clear}, consists in computing
\begin{equation}\label{refit:QO}
\hspace{-0.01cm}\,\tilde{x}^{\mathrm{QO}}\hspace{-0.15cm} \in\hspace{-0.1cm} \uargmin{x; \; \supp(\Gamma x) \subseteq \hat{\Ii}}\hspace{-0.12cm}
  \tfrac12 \| \Phi x \hspace{-0.03cm}-\hspace{-0.03cm} y \|^2 \hspace{-0.08cm}+\hspace{-0.1cm}
   \sum_{i\in\hat{\Ii}}\hspace{-0.1cm}\tfrac{\lambda}{2\norm{(\Gamma\hat x)_i}}\hspace{-0.05cm} \snorm{\hspace{-0.03cm}(\Gamma x)_i\hspace{-0.03cm}-\hspace{-0.03cm}\dotp{\hspace{-0.08cm}(\Gamma x)_i}{\hspace{-0.12cm}\tfrac{(\Gamma \hat x)_i}{\norm{(\Gamma \hat x)_i}}\hspace{-0.05cm}}\hspace{-0.1cm}\tfrac{(\Gamma \hat x)_i}{\norm{(\Gamma \hat x)_i}}\hspace{-0.05cm}}^2\hspace{-0.2cm}\enspace.
\end{equation}
This model promotes refitted solutions preserving to some extent the orientation
$\Gamma \hat x$ of the biased solution. It also shrinks the amplitude
of $\Gamma x$ all the more that the amplitude of $\Gamma \hat x$
are small.
As this model penalizes changes of orientation, we refer to it as QO for {\em Quadratic-penalized Orientation}.
This penalty does not promote any kind of direction preservation, and as for the HO model, contrast inversions may be observed between biased and refitted solutions (see Fig.~\ref{ex:TViso}(f)).
The quadratic term also over-penalizes large changes of orientation.

\section{Refitting with block penalties}
As mentioned in the previous section, the methods of  \cite{brinkmann2016bias} and \cite{deledalle2016clear}  have proposed
variational refitting formulations  that
not only aim at preserving the co-support $\hat{\Ii}$ but also
the orientation of $(\Gamma \hat x)_{i \in \hat{\Ii}}$.
In this paper, we propose to express these (two-steps) refitting procedures
in the following general framework
\begin{align}\label{general_refit}
  &
  \tilde{x}^{\phi} \in\uargmin{x; \; \supp(\Gamma x) \subseteq \hat{\Ii}}
  \tfrac12 \| \Phi x - y \|^2 +
  \sum_{i \in \hat{\Ii}} \phi((\Gamma x)_i, (\Gamma \hat x)_i)\enspace,
\end{align}
where $\phi : \RR^b \times \RR^b \to \RR$ is a refitting block penalty
($b \geq 1$ is the size of the blocks)
promoting $\Gamma x$ to share information with $\Gamma \hat{x}$ in some sense to be specified.
To refer to some features of the vector $\Gamma \hat{x}$, let us first define properly the notions of
relative orientation, direction and projection between two vectors.
\begin{defn}\label{def:cos_and_proj}
  Let $z$ and $\hat z$ being two vectors in $\RR^b$, we define
  \begin{align}   \cos(z, \hat{z})
  =
  \dotp{\tfrac{z}{\norm{z}}}{\tfrac{\hat z}{\norm{\hat z}}}=\tfrac{1}{\norm{z}\norm{\hat z}}\sum_{j=1}^b z_j\hat z_j \enspace,\\
  \qandq P_{\hat z}(z)= \dotp{z}{\tfrac{\hat{z}}{\norm{\hat{z}}}}\tfrac{\hat{z}}{\norm{\hat{z}}}=\tfrac{\norm{z}}{\norm{\hat{z}}} \cos(z, \hat{z})\hat{z}
  \enspace,
  \end{align}
  where $P_{\hat z}(z)$ is the orthogonal projection of $z$ onto $\Span(\hat z)$ (\ie the orientation axis of $\hat z$).
  We  say that
  $z$ and $\hat z$ share the same orientation (resp.~direction), if $|\cos(z, \hat{z})| = 1$ (resp.~$\cos(z, \hat{z}) = 1$).
\end{defn}
Thanks to Definition~\ref{def:cos_and_proj}, we can now reformulate the previous refitting models in terms of block penalties.
The Hard-constrained refitting models preserving Direction \eqref{refit_german:HD} and Orientation  \eqref{refit_german:HO} of  \cite{brinkmann2016bias} as well as the flexible Quadratic  model of CLEAR \cite{deledalle2016clear} correspond to the following block penalties:
\begin{align}\label{pen:HO}
 & \phi_{\mathrm{HD}}(z, \hat{z})
  ={\iota_{\{z \in \RR^b : \cos(z, \hat{z}) = 1\}}}
  \\
  \label{pen:HD}
&  \phi^{}_{\mathrm{HO}}(z, \hat{z})
  ={\iota_{\{z \in \RR^b : |\cos(z, \hat{z})|= 1\}}}
  \\
  \label{pen:QO}
    \qandq &
  \phi_{\mathrm{QO}}(z, \hat z)
  =
  \tfrac{\lambda}{2\norm{\hat{z}}}
 \left|\left|z-P_{\hat z}(z)\right|\right|^2
  =
  \tfrac{\lambda}{2}
  \tfrac{\norm{z}^2}{\norm{\hat{z}}}(1-
  \cos^2(z, \hat{z}))\enspace.
\end{align}
where $\iota_{\mathcal{C}}$ is the $0/+\infty$ indicator function of a set $\mathcal{C}$.
These block penalties are either insensitive to directions (HO and QO) or intolerant to small changes of orientations (HD and HO), hence not satisfying (\lcf drawbacks visible in Fig.~\ref{ex:TViso}).

When $b=1$, the orientation-based penalties (QO and HO) have absolutely no effect while the direction-based penalty HD preserves the sign of $(\Gamma \hat x)_i$.
In this paper, we argue that the direction of $(\Gamma \hat x)_i$, for any $b\geq 1$, carries important information that is worth preserving when refitting,
at least to some extent.

\subsection{Desired properties of refitting block penalties}
To compute global optimum of the refitting model \eqref{general_refit}, we only consider convex refitting block penalties $z \mapsto \phi(z, \hat{z})$.
Hence, we now introduce properties a block penalty $\phi$ should satisfy for refitting purposes:
\begin{itemize}\setlength{\itemindent}{.1in}
\item [(P1)] $\phi$ is convex, non negative and $\phi(z, \hat{z}) = 0$, if $\cos(z, \hat z) = 1$ or $\norm{z}=0$,
\item [(P2)] $\phi(z', \hat{z}) \geq \phi(z, \hat{z})$
  if $\norm{z'} = \norm{z}$ and $\cos (z, \hat{z}) \geq \cos (z', \hat{z})$,
\item [(P3)] $z \mapsto \phi(z, \hat{z})$ is continuous,
\end{itemize}
Property (P1) stipulates that no configuration can be more favorable than $z$ and $\hat z$ having the same direction. 
Hence, the direction of the refitted solution should be encouraged to follow the one of the biased solution.
Property (P2) imposes that for a fixed amplitude, the penalty should be increasing w.r.t. the angle formed with $\hat z$.
Property (P3) enforces refitting that can continuously adapt to the data and be robust to small perturbations.\vspace{-0.em}

\begin{table}[!ht]
  \centering
  \caption{\label{tab:properties}Properties satisfied by the considered block penalties $\phi$.}
  \begin{tabular}{c@{\hspace{.5cm}}c@{\hspace{.5cm}}c@{\hspace{.5cm}}c@{\hspace{.5cm}}cc}
    \hline
    Properties & HO & HD & QO& SD
    \\
    \hline
    \hline
    1 & $\surd$ & $\surd$ & $\surd$  &   $\surd$ 
    \\
    2 &  &  $\surd$ & & $\surd$ 
    \\
    3 &        &          & $\surd$ & $\surd$
    \\
    \hline

  \end{tabular}
  \vspace*{-0.em}
\end{table}

\subsection{A new flexible refitting block penalty}

We now introduce our refitting block penalty
designed to preserve the desired features of $\hat{z}=\Gamma \hat x$ in a simple way.
The  {\em Soft-penalized Direction} penalty reads
\begin{align}\label{pen:SD}
\phi^{}_{\mathrm{SD}}(z, \hat{z}) =
\lambda \norm{z}(1 - \cos (z, \hat{z}))\enspace.
\end{align}
The properties of the different studied block penalties are presented in Table \ref{tab:properties}. The proposed SD model is the only one satisfying  all the desired properties.
As illustrated in Fig.~\ref{fig:BP}, it is a continuous penalization that increases continuously with respect to the absolute angle between $z$  and $\hat z$.

\newpage

\newcommand{\scaleX}{0.255\linewidth}
\newcommand{\scaleXB}{0.25\linewidth}
\newcommand{\scaleY}{0.27\linewidth}
\newcommand{\scaleZ}{0.275\linewidth}
\newcommand{\sidecapY}[1]{{\begin{sideways}\parbox{\scaleY}{\centering #1}\end{sideways}}}
\begin{figure}[!t]
  \begin{center}%
    \hfill%
\sidecapY{QO}\hfill\hspace{0.00\linewidth}%
{\includegraphics[width=\scaleX]{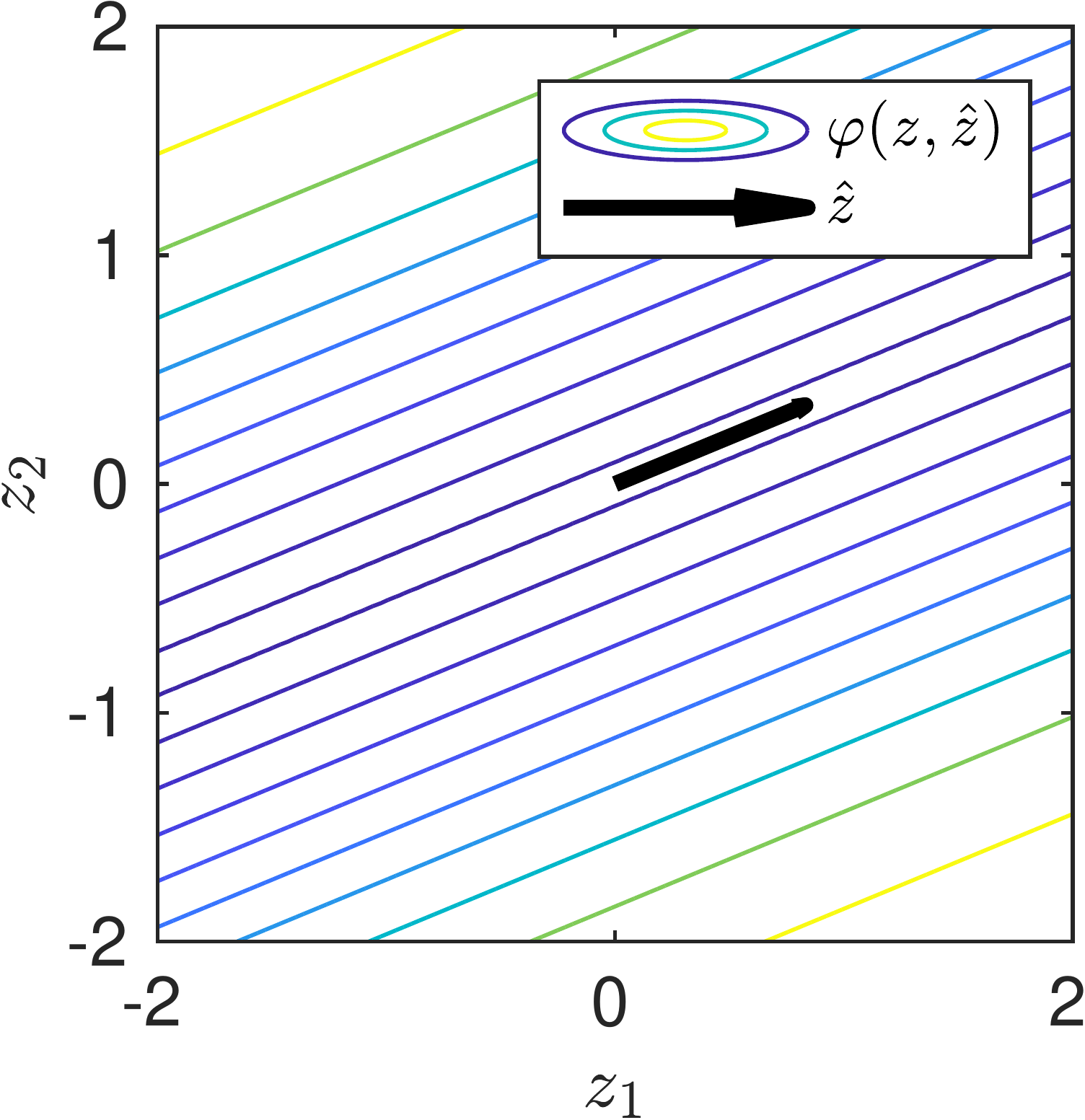}}\hfill%
{\includegraphics[width=\scaleY]{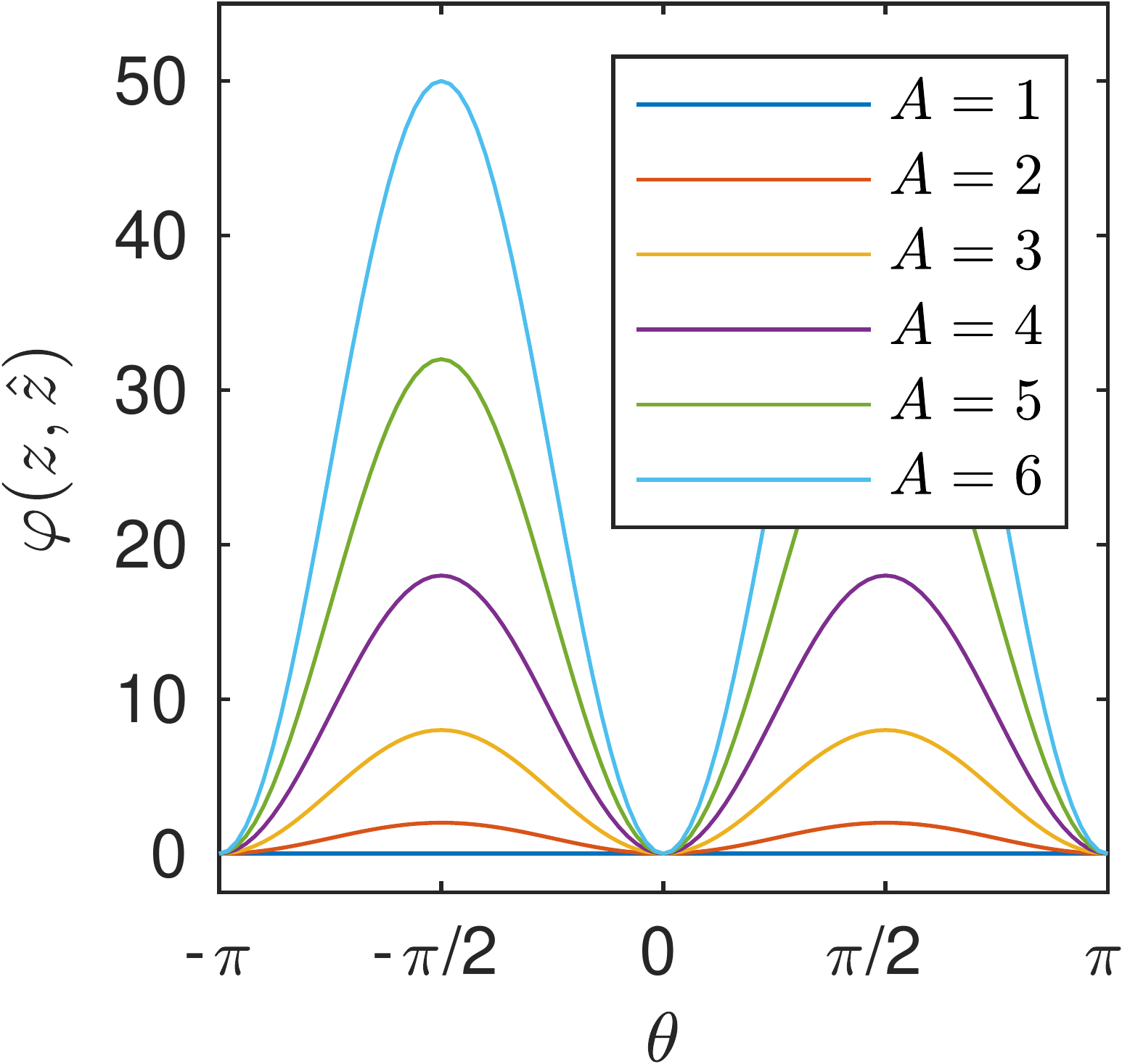}}\hfill%
{\includegraphics[width=\scaleZ]{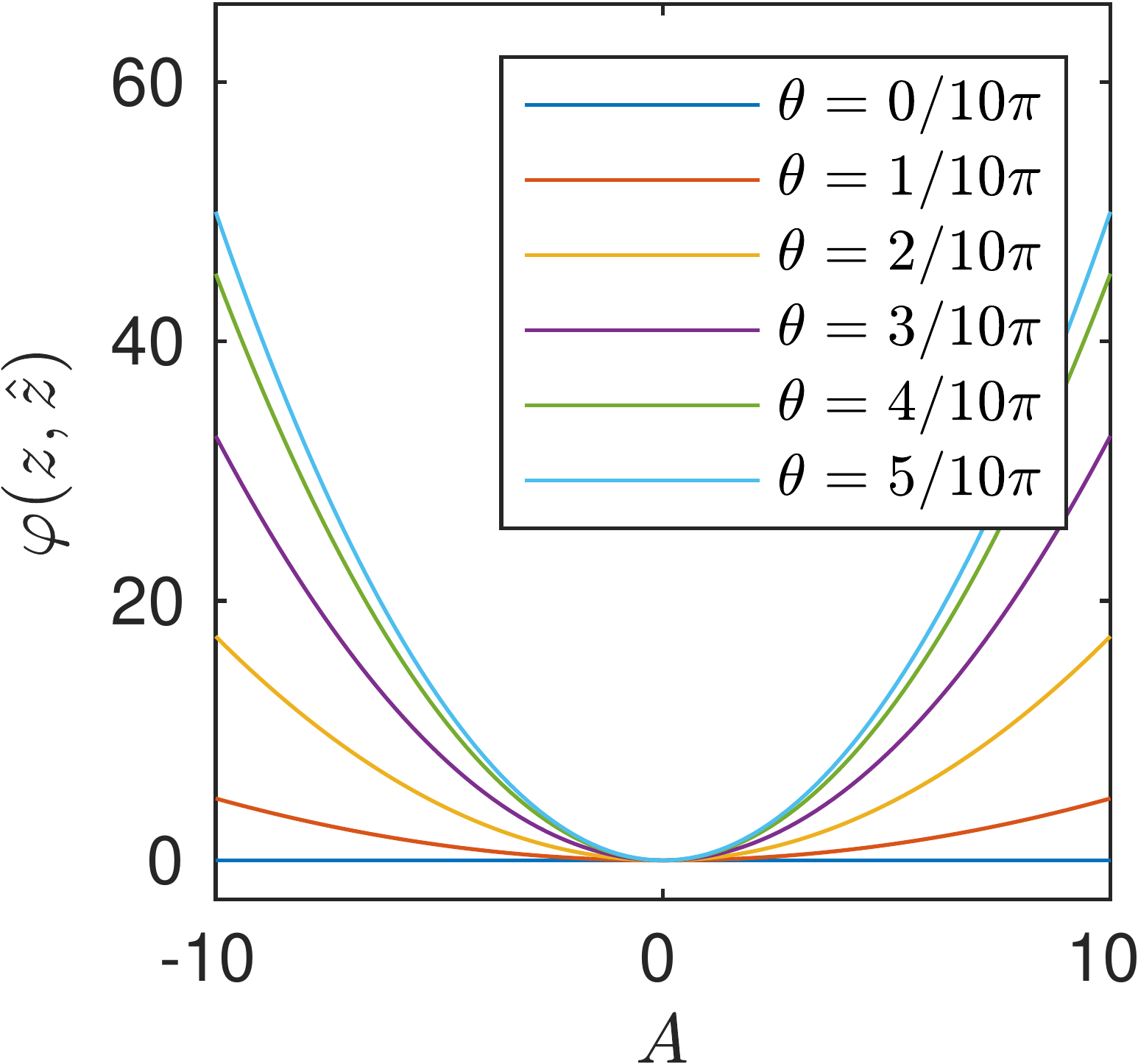}}\hfill%
    { }
    \\%
    \hfill%
\sidecapY{SD}\hfill\hspace{0.005\linewidth}%
{\includegraphics[width=\scaleXB]{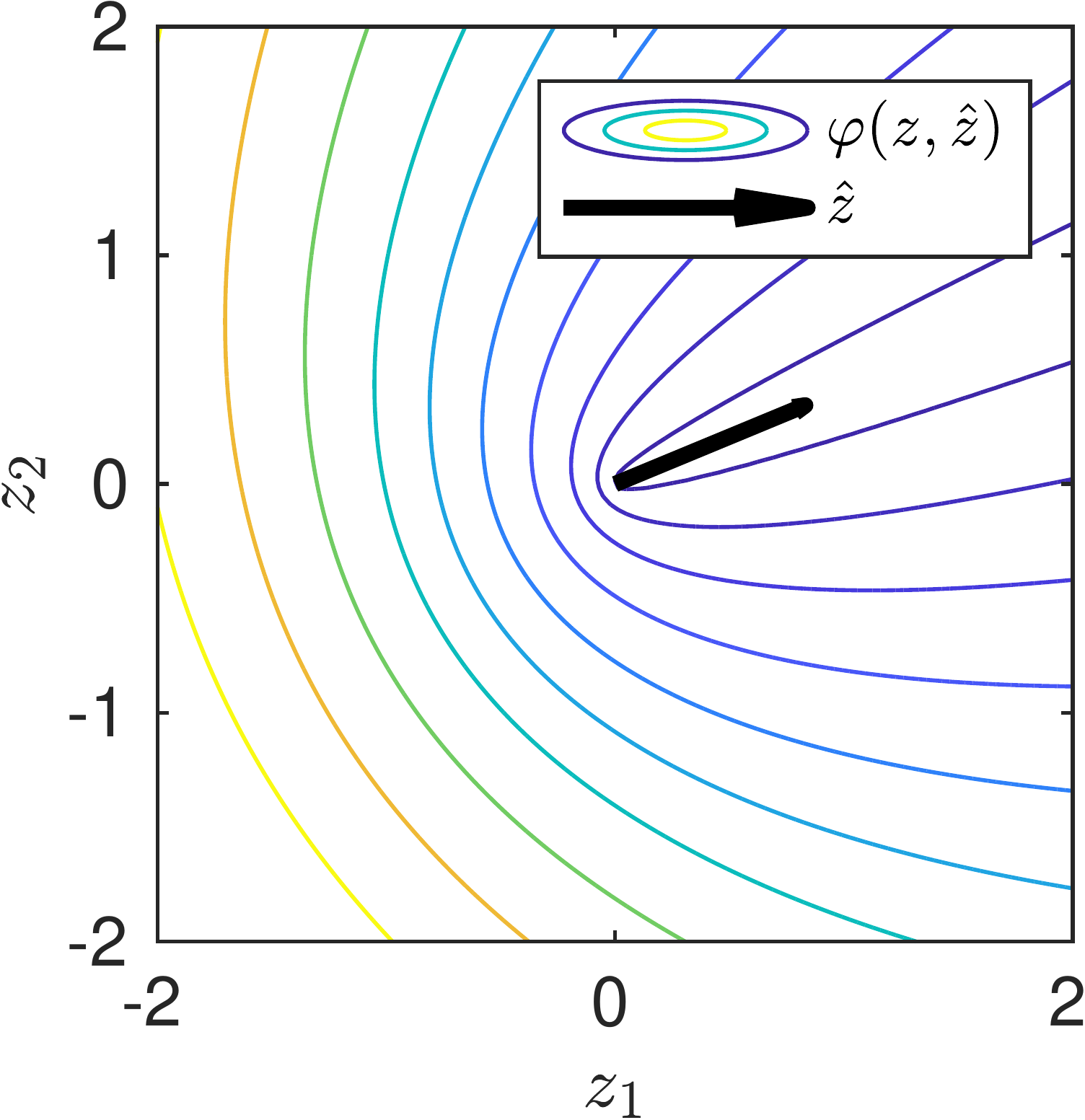}}\hfill%
{\includegraphics[width=\scaleY]{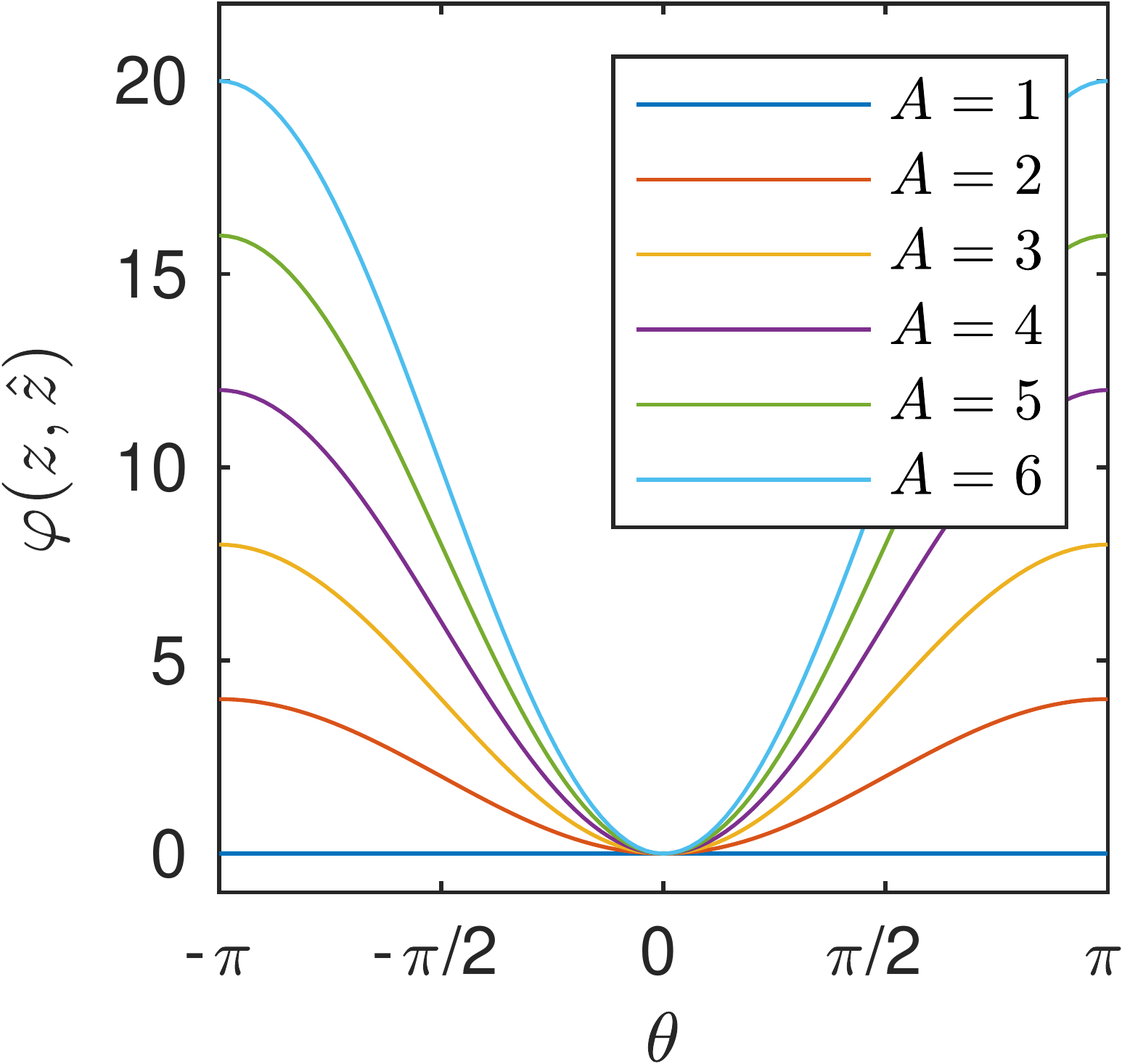}}\hfill%
{\includegraphics[width=\scaleZ]{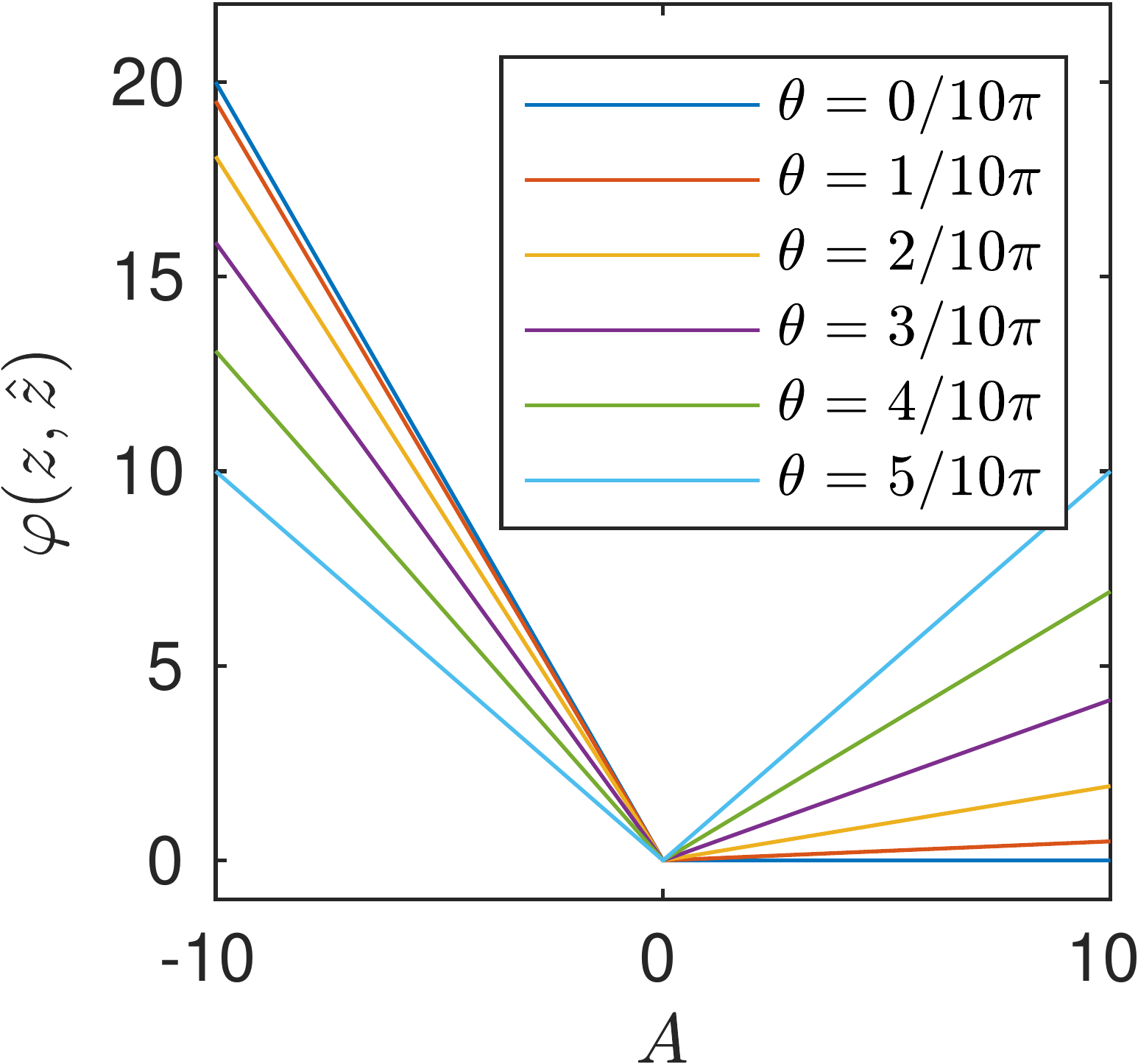}}\hfill%
    { }
    \\
  \end{center}
  \vspace{-2em}

  \newsavebox{\smlmat}
  \savebox{\smlmat}{$\left(\begin{smallmatrix} \cos \theta & -\sin \theta \\ \sin \theta & \cos \theta \end{smallmatrix}\right)$}
  \caption{\label{fig:BP} Illustration of continuous block penalties QO and SD:
    (left) 2D level lines of  $\phi$ for
    $z=(z_1,z_2) = A $\usebox{\smlmat}$\hat{z}$,
    (middle) evolution regarding $\theta$ and (right) $A$.}

\end{figure}

\subsection{SD block penalty: the best of both Bregman worlds}
Denoting $D_i(\Gamma x)=D_{ \norm{\cdot}} ((\Gamma x)_i,(\Gamma \hat x)_i)$
as introduced in Section~\ref{sec:intro_breg}, the refitting models  given in  \eqref{iterative_refit} and \eqref{refit_german:HD} can be expressed as\vspace{-0.1cm}
\begin{align}
\label{iterative_refit2}
\tilde{x}^{\mathrm{IB}}&\in\uargmin{x} \tfrac12 \|\Phi x - y \|^2 +\lambda\sum_{i=1}^m D_i(\Gamma x)\enspace,\\
\label{refit_german:HD3}
  \tilde{x}^{\mathrm{HD}} &\in \uargmin{x}
  \tfrac12 \| \Phi x - y \|^2 &\mathrm{s.t. }\,\,  D_i(\Gamma x)=0,\,\forall i\in [m]\enspace.
\end{align}

Observing that the new SD block penalty \eqref{pen:SD} may be rewritten  as
$
      \phi^{}_{\mathrm{SD}}(z, \hat{z})
  =
\lambda
  \left(
  \norm{z} - \dotp{
    z
  }{
   {\hat{z}}/{\norm{\hat{z}}}
  }
  \right)=\lambda D_{\norm{\cdot}}(z,\hat z),
$
 the SD refitting model is
\begin{align}\label{SD:bregman}
  \tilde{x}^{\mathrm{SD}}\in\uargmin{x}
  \tfrac12 \| \Phi x - y \|^2 +\lambda \sum_{i\in\hat{\Ii}} D_i(\Gamma x),\hspace{0.3cm} &\mathrm{s.t. }\hspace{0.1cm} D_i(\Gamma x)=0,\,\forall i\in\hat{\Ii}^c\enspace.
\end{align}
With such reformulations, connections between refitting models \eqref{iterative_refit2}, \eqref{refit_german:HD3} and \eqref{SD:bregman} can be clarified.
The solution $\tilde{x}^{\mathrm{IB}}$ \cite{Osher} is too relaxed, as it only penalizes the directions $(\Gamma  x)_i$ using $(\Gamma \hat x)_i$, without aiming at preserving the co-support of $\hat x$.
The solution $\tilde{x}^{\mathrm{HD}}$ \cite{brinkmann2016bias} is too constrained: the direction within the co-support is required to be preserved exactly. Our proposed refitting $\tilde{x}^{\mathrm{SD}}$ lies in-between: it preserves the co-support, while authorizing some directional flexibility, as illustrated by the sharper square edges in Figure 1(g).

With respect to the Hard-constrained approach \cite{brinkmann2016bias}, an important difference is that we consider local inclusions of subgradients of the function $\lambda\norm{\cdot}_{1,2}$ at point $(\Gamma x)_i$ instead of the global inclusion of subgradients of the function $\lambda\norm{\Gamma \cdot}_{1,2}$ at point $x$ as in HD \eqref{refit_german:HD} and HO \eqref{refit_german:HO}. Such a change of paradigm allows to adapt the refitting locally by preserving the co-support while including the flexibility of the original Bregman approach \cite{Osher}.

\section{Refitting in practice}
We first describe how computing $\hat x$ and its refitting $\tilde x$ in two successive steps and then propose a new numerical scheme for the joint computation of $\hat x$ and  $\tilde x$.\vspace{-0.1cm}
\subsection{Biased problem and posterior refitting}
To obtain $\hat x$ solution of \eqref{isotv},  we  consider the primal dual formulation that reads
\begin{align}\label{isotv:pd}
\min_{x}\max_z \tfrac12 \norm{\Phi x - y }^2 + \langle \Gamma x,z\rangle - \iota_{B_2^\lambda}(z)\enspace,
\end{align}
where $\iota_{B_2^\lambda}$ is the indicator function of the $\ell_2$ ball of radius $\lambda$ that is $0$ if $\norm{z_i}\leq \lambda$ for all $i\in[m]$ and $+\infty$ otherwise.
This problem can be solved  with the iterative primal-dual algorithm \cite{CP} presented in the left part of  Algorithm \eqref{algo_final}. For positive parameters satisfying $\tau\sigma<\|\Gamma\|$ and $\theta\in[0,1]$,
the iterates $(\hat z^k,\hat x^k)$ converge to
a saddle point $(\hat z,\hat x)$ satisfying $\hat z_i\in\partial \norm{\cdot}((\Gamma \hat x)_i)$, $\forall i$.

Assume that the  co-support $\hat{\Ii} = \supp(\Gamma\hat{x})$  of the biased solution has been identified,   a posterior refitting can be obtained by solving \eqref{general_refit} for  any refitting block penalty $\phi$.
To that end, we write the characteristic function of  co-support preservation as
$\sum_{i\in\hat{\Ii}^c} \iota_{\{0\}}(z_i),$
where $\iota_{\{0\}}(z)=0$ if  $z=0$ and $+\infty$ otherwise. By introducing the convex function
\begin{equation}\label{def:omega}
\omega_{\phi}(z,\hat z,\hat{\Ii} )= \sum_{i\in\hat{\Ii}^c} \iota_{\{0\}}(z_i)+\sum_{i \in \hat{\Ii}} \phi(z_i, \hat z_i)\enspace,\vspace{-0.1cm}
\end{equation}
the general refitting problem \eqref{general_refit} can be expressed as
\begin{equation}\label{general_refit_reform}
  \tilde{x}^{\phi} \in\uargmin{x}
  \tfrac12 \| \Phi x - y \|^2 + \omega_{\phi} (\Gamma x,\Gamma \hat x,\hat{\Ii})\enspace.\vspace{-0.1cm}
\end{equation}
Subsequently, we can consider its primal dual formulation
\begin{equation}\label{general_refit_reform:pd}
  \min_{x}\max_{z}
  \tfrac12 \| \Phi x - y \|^2 +\langle \Gamma x,z\rangle -\omega^*_{\phi} (z,\Gamma \hat x,\hat{\Ii}) \enspace,\vspace{-0.1cm}
\end{equation}
where $\omega_\phi^*(z,\Gamma \hat x,\hat{\Ii})=\sup_{x} \langle z,x\rangle-\omega_{\phi}(x,\Gamma \hat x,\hat{\Ii})$ is the convex conjugate, with respect to the first argument, of $\omega_\phi(\cdot,\Gamma \hat x,\hat{\Ii})$.
Such problem can again be solved with the primal-dual algorithm \cite{CP}.
The crucial point is to have an accurate estimation of the vector $\Gamma \hat x$ and its support $\hat{\Ii}$. Yet, it is well known that estimating $\supp(\Gamma \hat x)$ from an estimation $\hat x^k$ is not stable numerically:
the support $\supp(\Gamma \hat x^k)$ can be far from $\supp(\Gamma \hat x)$  even though $\hat x^k$ is arbitrarily close to $\hat x$.

\subsection{Joint-refitting algorithm}

We now  introduce a general algorithm aiming to jointly solve the original problem \eqref{isotv} and the refitting one \eqref{general_refit} for any refitting block penalty $\phi$. 
This framework has been developed for stable projection onto the support in \cite{Deledalle_Papadakis_Salmon15} and later extended to refitting with the Quadratic Orientation penalty in \cite{deledalle2016clear}.
The strategy consists in solving in parallel the two problems \eqref{isotv:pd} and \eqref{general_refit_reform:pd}.

Two  iterative primal-dual algorithms are used for the biased variables $(\hat z^k,\hat x^k)$ and the refitted ones $(\tilde z^k,\tilde x^k)$. 
Let us now present the whole algorithm:
\begin{equation}\label{algo_final}\normalsize
  \begin{array}{@{}ll@{\hspace{.4em}}|@{\hspace{.4em}}ll@{}}
    \\[-1.4em]
    \hat z^{k+1}_i&=\frac{\hat z^k_i+\sigma( \Gamma \hat v^k)_i}{\max(\lambda,\norm{\hat z^k_i+\sigma (\Gamma \hat v^k)_i})}&
   \, \hat{\Ii}^k&=\enscond{i \in [m]}{\norm {\hat z^k_i+\sigma (\Gamma \hat v^k )_i}>\lambda}\\
 &&\, \zt^{k+1}&=\prox_{\sigma \omega_\phi^*(\cdot,\Psi(\hat z^k,\hat v^k),\hat{\Ii}^k )}(\zt^k+\sigma \Gamma \tilde v^k)\\
  \hat \x^{k+1}&=\Phi_\tau ^{-}\left(\hat\x^k+\tau( \Phi^\top y-\Gamma^\top \hat z^{k+1} )\right)
 &\, \xt^{k+1}&=\Phi_\tau ^{-}\left(\xt^k+\tau(\Phi^\top y-\Gamma^\top\zt^{k+1})\right),\\
    \hat v^{k+1}&=\hat\x^{k+1}+\theta(\hat\x^{k+1}-\hat\x^k)&
  \,\tilde v^{k+1}&=\xt^{k+1}+\theta(\xt^{k+1}-\xt^{k}),
\vspace{-1em}
\end{array}
\vspace{1em}
\end{equation}
with the operator $\Phi_\tau ^{-}=(\Id+\tau\Phi^\top\Phi)^{-1}$ and the auxiliary variables $\hat v^k$ and $\tilde v^k$ that  accelerate the algorithm.
Following \cite{CP}, for  any positive scalars $\tau$ and  $\sigma$ satisfying $\tau\sigma \|\Gamma^\top\Gamma\|<1$ and $\theta\in[0,1]$, the estimates $(\hat z^k,\hat x^k,\hat v^k)$ of the biased solution converge to $(\hat z, \hat x,\hat x)$, where $(\hat z,\hat x)$ is  a saddle point  of \eqref{isotv:pd}.

In the right part of  Algorithm \eqref{algo_final},  we rely on
 the proximal operator  that is, for a convex function $\psi$ and at point $z_0$,
$\prox_{\sigma \psi}(z)=\argmin_{z}\tfrac1{2\sigma}\norm{z-z_0}^2+\psi(z)$.
%
From the block structure of the function $\omega_\phi$ defined in \eqref{def:omega},  the computation of its  proximal operator  may be realized  pointwise. Since $\iota_{\{0\}}(z)^*=0$, we have
\begin{align}
\prox_{\sigma \omega_\phi^*}(z^0,\Gamma \hat x,\hat{\Ii})_i&=
 \begin{cases}
  \prox_{\sigma \phi^*}(z^0_i,(\Gamma \hat x)_i), &\quad \textrm{if } i\in \hat{\Ii}\enspace,\\
    z^0_i,&\quad \textrm{otherwise}\enspace.
 \end{cases}
\end{align}
Table \ref{tab:prox_phi_dual} gives the expressions of the dual functions $\phi^*$ with respect to their first variable and their related proximal operators $\prox_{\sigma \phi^*}$ for the refitting block penalties considered in this paper. All details are given in the Appendix.

\begin{table}[!t]
  \centering
  \caption{\label{tab:prox_phi_dual} Convex conjugates and proximal operators of the studied block penalties $\phi$. 
  }
  \begin{tabular}{l@{\hspace{0.02cm}}|l@{\hspace{0.02cm}}|l}
    \hline
    \multicolumn{1}{c|}{$\phi$}& \multicolumn{1}{c|}{$\phi^*(z,\hat z)$}& \multicolumn{1}{c}{$\prox_{\sigma \phi^*}(z_0,\hat z)$}\\
    \hline
    \hline
    HO&$
    \choice{
      0,
      & \quad\quad \ifq\cos (z, \hat{z}) =0\\
      +\infty,
      & \quad\quad
      \otherwise
    }$
    &$z_0 -P_{\hat z}(z_0)$\\
    \hline
    HD&
    $\choice{
      0,
      & \quad\quad \ifq \cos(z, \hat{z}) \leq 0\\
      +\infty,
      & \quad\quad
      \otherwise
    }$
    &$ \choice{
      z_0 -P_{\hat z}(z_0),
      & \quad \ifq \langle z_0,\hat z\rangle\geq 0\\
      z_0,
      &\quad
      \otherwise
    }$\\
    \hline
    QO&$\choice{
      \tfrac{\norm{\hat z}}{2\lambda}\norm{z}^2, &\;\;\, \ifq \cos(z, \hat{z})  =0\\
      +\infty,
      &\;\;\,
      \otherwise
    }$&$\tfrac{\lambda}{\lambda+\sigma \norm{\hat z}}\left(z_0 -P_{\hat z}(z_0)\right)$\\
    \hline
    SD&
    $ \choice{
      0,
      & \quad\quad \ifq\hspace{-0.2cm}\norm{z+\lambda \tfrac{\hat z}{\norm{\hat z}}}\leq \lambda {}^\dag\\
      +\infty,
      &\quad\quad
      \otherwise
    }$
    &$\lambda\left(\tfrac{z_0+\lambda\tfrac{\hat z}{\norm{\hat z}}}{\max(\lambda,\norm{z_0+\lambda \tfrac{\hat z}{\norm{\hat z}}})}-\tfrac{\hat z}{\norm{\hat z}}\right)$\\
    \hline
  \end{tabular}
  \\
    {\scriptsize ${}^\dag$: note that the condition implies that $\cos(z,\hat z)\leq 0$.}
    \vspace{-1em}
\end{table}

The idea behind  this joint-refitting algorithm is to perform online co-support detection using  the dual variable $\hat z^k$ of the biased variable $\hat x^k$.
 From  relations in \eqref{opt_z}, we expect  at convergence $\hat z^k$  to   saturate  on the  support of $\Gamma \hat x$ and to satisfy the optimality condition
$\hat z_i^k=\lambda\tfrac{(\Gamma \hat x)_i}{\norm{(\Gamma \hat x)_i}}$.
 In practice, the norm of the dual variable $\hat z^k_i$ saturates to $\lambda$ relatively fast onto $\hat{\Ii}$.
As a consequence,  it is far more stable to detect the support of $\Gamma \hat x$ with the dual variable $\hat z^k$ than with the vector $\Gamma \hat x^k$ itself.
In the first step of Algorithm \eqref{algo_final}, the condition $\norm{\hat z^k_i+\sigma (\Gamma \hat v^k )_i}>\lambda$ is thus  used to detect elements of the support $i\in\hat{\Ii}^k$ of $\Gamma\hat x^k$ along iterations\footnote{As  in \cite{brinkmann2016bias}, extended support  $\norm{\hat z_i}=\lambda$ can be tackled  by testing  $\norm {(\hat z^k+\sigma \Gamma \hat v^k )_i}\geq\lambda$.}.
 The function
  $\Psi(\hat z^k,\hat v^k)$ aims at approximating $\Gamma \hat x$ with the current values of the available variables $(\hat z^k,\hat v^k)$ of the biased problem that is solved simultaneously.
Following \cite{deledalle2016clear}, the function $\Psi$ can be chosen as
\begin{align}
  \Psi(\hat z^k, \hat v^k)_i
  =
  \tfrac{\norm{\hat{\nu}_i^k} - \lambda}{\sigma \norm{\hat{\nu}_i^k}}
  \hat{\nu}_i^k
  \qwhereq
  \hat{\nu}_i^k = (\hat z^k+\sigma \Gamma \hat v^k)_i\enspace.
\end{align}
that  satisfies $(\Psi(\hat z,\hat v))_i=(\Gamma \hat x)_i$ at convergence, while appearing  to give very stable online estimations of  the direction of $\Gamma \hat{x}$ through $\hat z^k$.

This joint-estimation  considers at every iteration $k$  different  refitting functions $\omega_\phi^*(.,\Psi(\hat z^k, \hat v^k),\hat{\Ii}^k)$ in  \eqref{general_refit_reform}. For fixed values of $\Psi(\hat z^k, \hat v^k)$ and $\hat{\Ii}^k$, the refitted variables $(\tilde z^k,\tilde x^k)$ in the right part of the Algorithm \eqref{algo_final}
converges since it exactly corresponds to the primal-dual algorithm \cite{CP} applied to the problem \eqref{general_refit_reform:pd}.
However, unless $b=1$ (see \cite{Deledalle_Papadakis_Salmon15}),
we do not have guarantee of convergence of the presented scheme with a varying $\omega_\phi^*$.
As in \cite{deledalle2016clear}, we nevertheless observe convergence and
a very stable behavior for this algorithm.

In addition to its better numerical stability,  the running time of  joint-refitting  is more interesting than the posterior approach. 
In Algorithm \eqref{algo_final},  the refitting variables at iteration $k$ require the biased variables at the same iteration and the whole process  can be realized in parallel without significantly affecting the running time of the  original biased process. On the other hand, posterior refitting is necessarily sequential and  the running time is doubled in general.

\section{Results}

\begin{figure}[!t]
  \subfigure[Noisy {\scriptsize (22.10)}]{%
    \begin{minipage}{.193\linewidth}
    \includegraphics[width=1\linewidth]{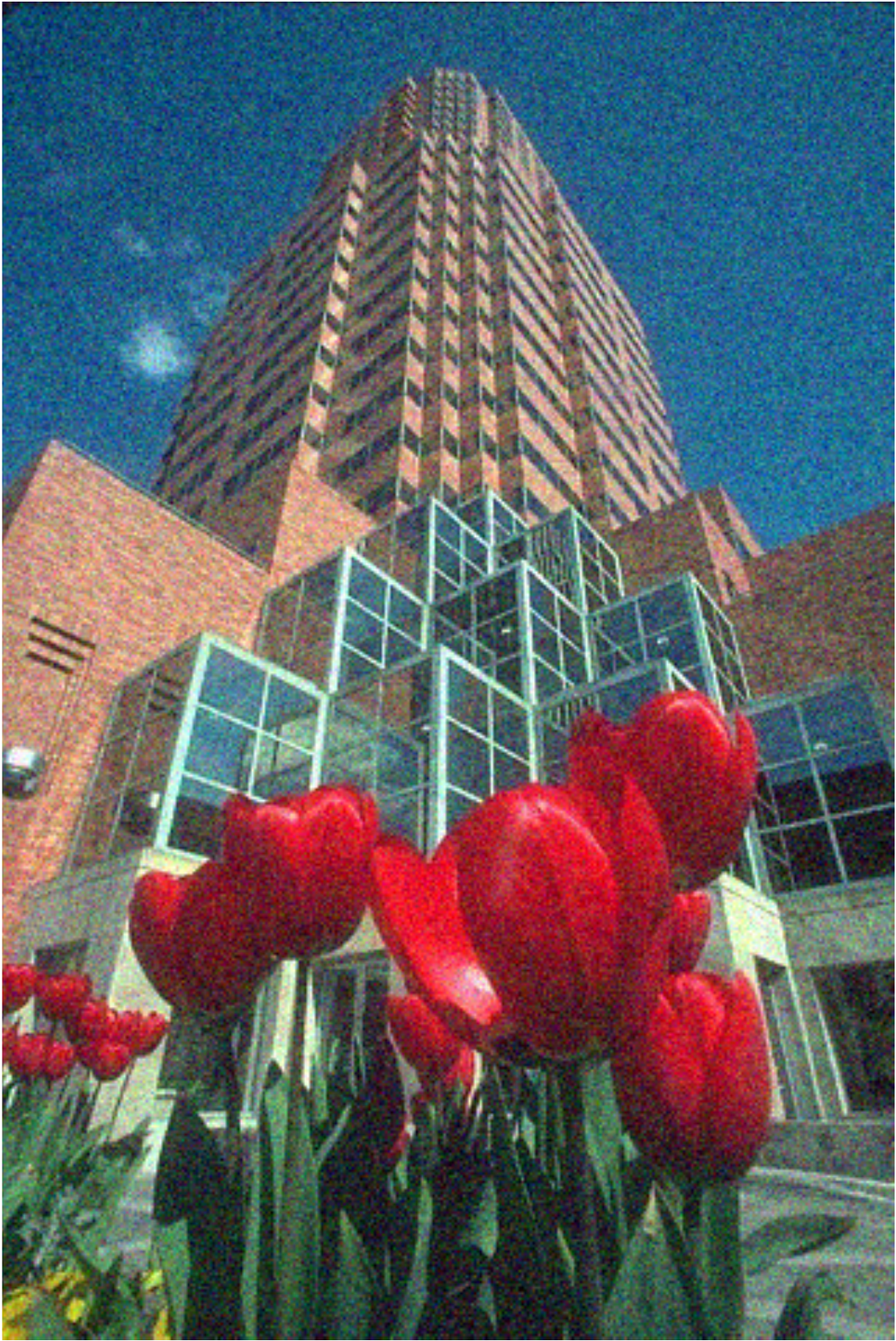}\\[.05em]
      \includegraphics[width=.49\linewidth,viewport=0 140 85 225,clip]{f_subimg1}\hfill%
      \includegraphics[width=.49\linewidth,viewport=135 370 210 445,clip]{f_subimg1}\\[-.5em]
    \end{minipage}%
  }\hfill
  \subfigure[TViso {\scriptsize (23.28)}]{%
    \begin{minipage}{.193\linewidth}
    \includegraphics[width=1\linewidth]{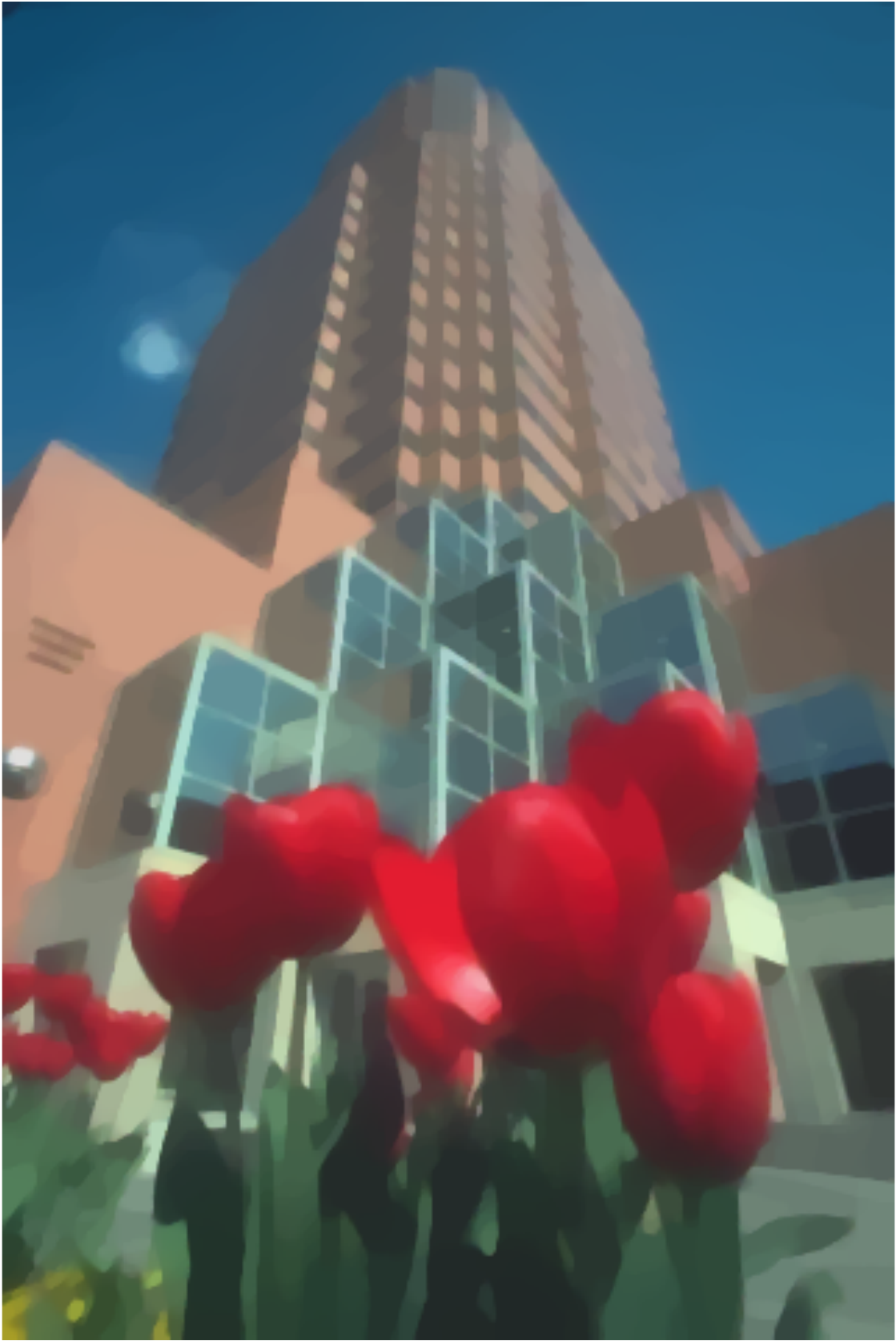}\\[.05em]
      \includegraphics[width=.49\linewidth,viewport=0 140 85 225,clip]{f_subimg2}\hfill%
      \includegraphics[width=.49\linewidth,viewport=135 370 210 445,clip]{f_subimg2}\\[-.5em]
    \end{minipage}%
  }\hfill
  \subfigure[HD {\scriptsize (23.75)}]{%
    \begin{minipage}{.193\linewidth}
    \includegraphics[width=1\linewidth]{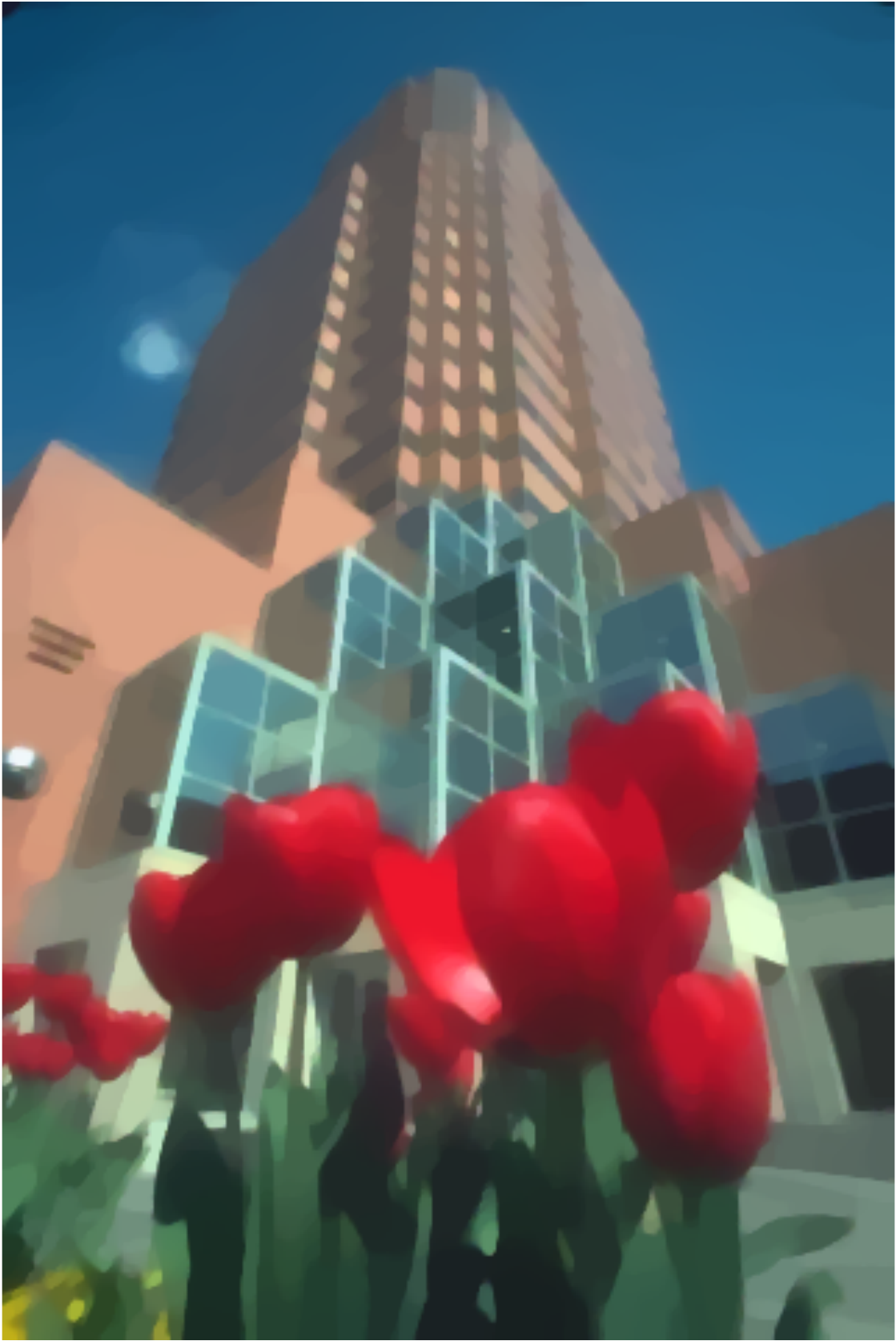}\\[.05em]
      \includegraphics[width=.49\linewidth,viewport=0 140 85 225,clip]{f_subimg4}\hfill%
      \includegraphics[width=.49\linewidth,viewport=135 370 210 445,clip]{f_subimg4}\\[-.5em]
    \end{minipage}%
  }\hfill
  \subfigure[QO {\scriptsize (26.12)}]{%
    \begin{minipage}{.193\linewidth}
      \includegraphics[width=1\linewidth]{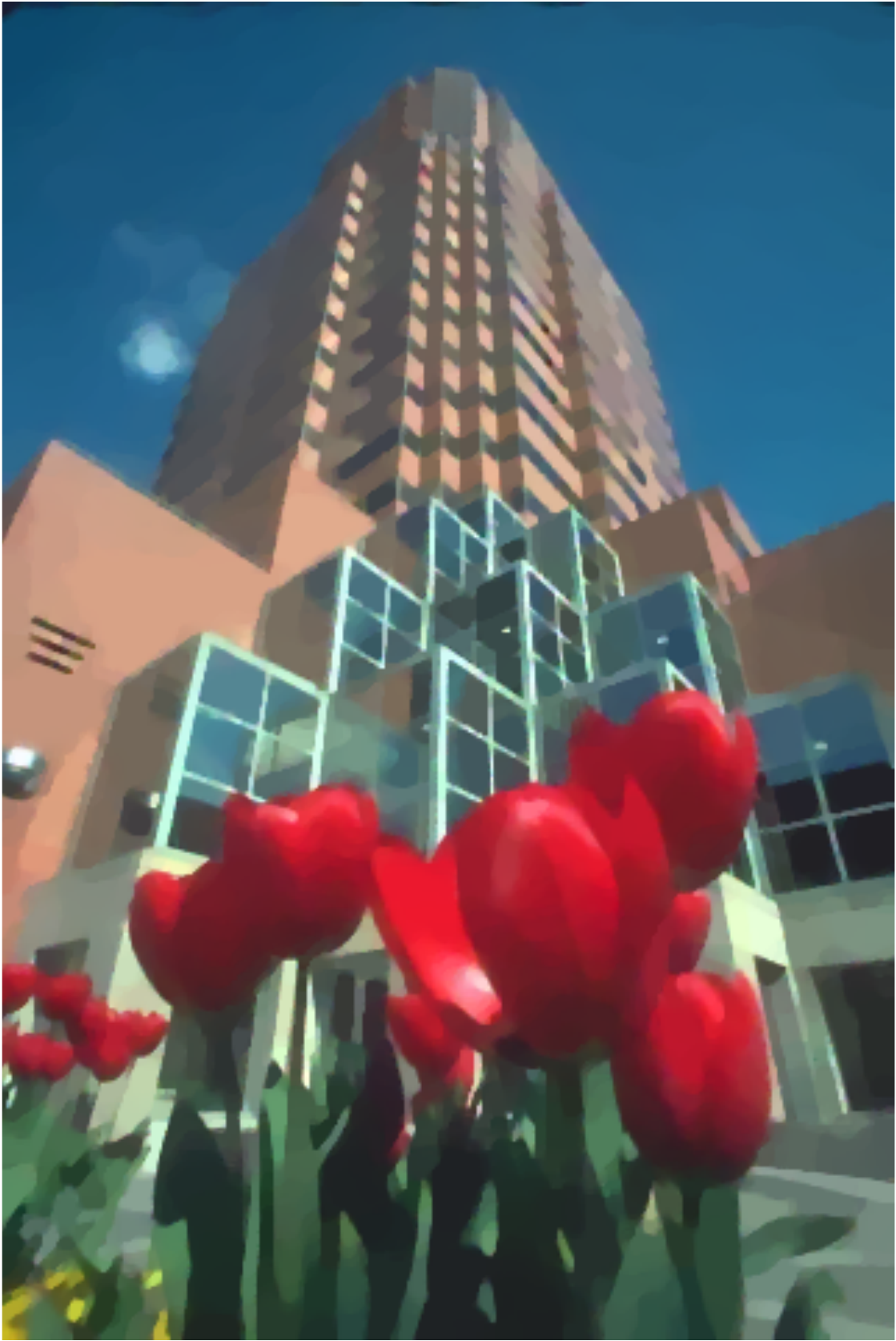}\\[.05em]
      \includegraphics[width=.49\linewidth,viewport=0 140 85 225,clip]{f_subimg5}\hfill%
      \includegraphics[width=.49\linewidth,viewport=135 370 210 445,clip]{f_subimg5}\\[-.5em]
    \end{minipage}%
  }\hfill
  \subfigure[SD~{\scriptsize (27.68)}]{%
    \begin{minipage}{.193\linewidth}
    \includegraphics[width=1\linewidth]{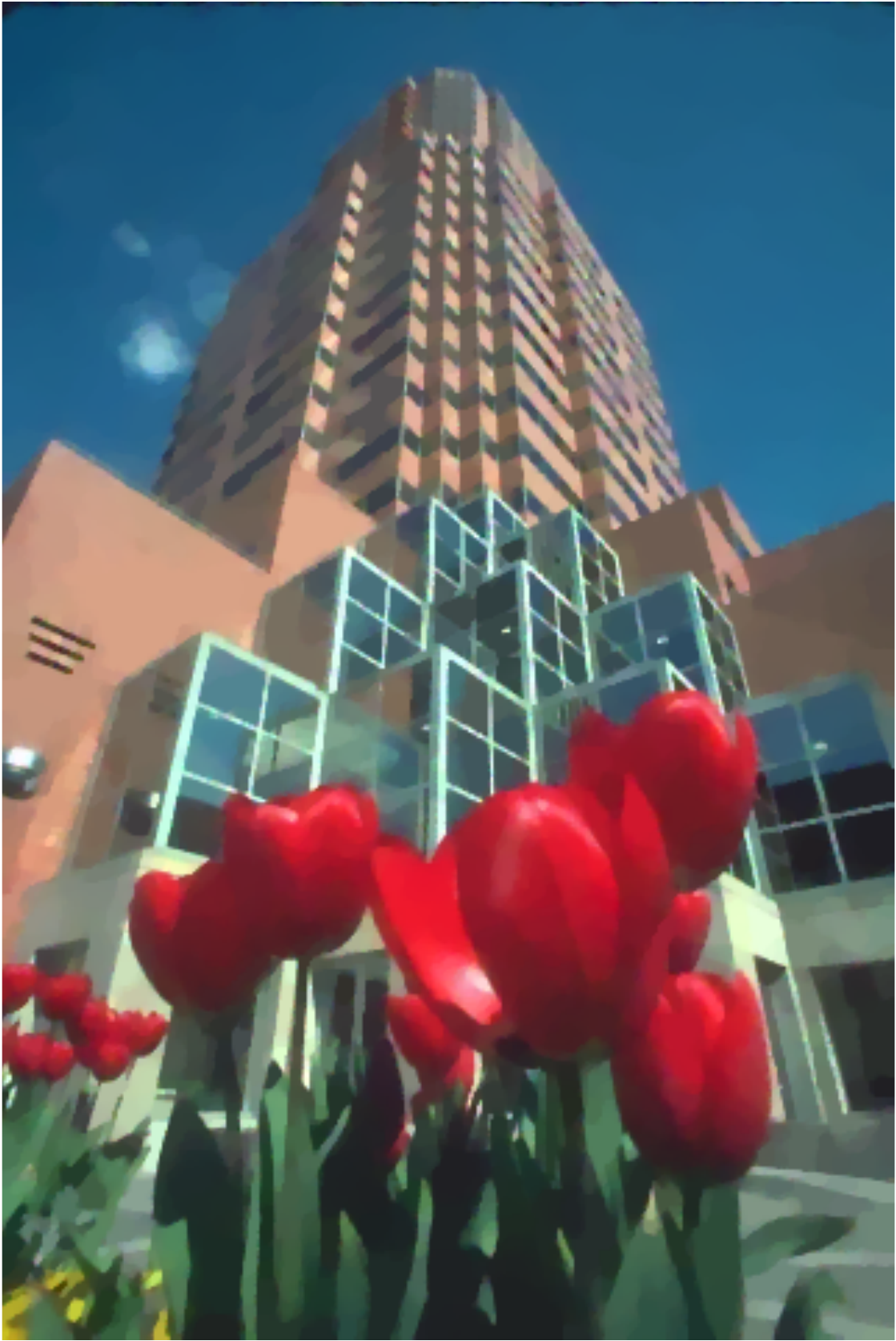}\\[.05em]
    \includegraphics[width=.49\linewidth,viewport=0 140 85 225,clip]{f_subimg6}\hfill%
    \includegraphics[width=.49\linewidth,viewport=135 370 210 445,clip]{f_subimg6}\\[-.5em]
    \end{minipage}%
  }\\
  \vspace{-1em}
  \caption{
    (a) An 8bit color image corrupted by Gaussian noise with standard deviation $\sigma=20$.
    (b) Solution of TViso.
    Debiased solution with (c) HD, (d) QO and (e) SD.
    The Peak Signal to Noise Ratio (PSNR) is indicated in brackets bellow each image.}
  \label{fig:denoising}
  \vspace{-1em}
\end{figure}

\begin{figure}[!t]
  \subfigure[Blurry {\scriptsize (23.14)}]{%
    \begin{minipage}{.470\linewidth}
    \includegraphics[width=1\linewidth]{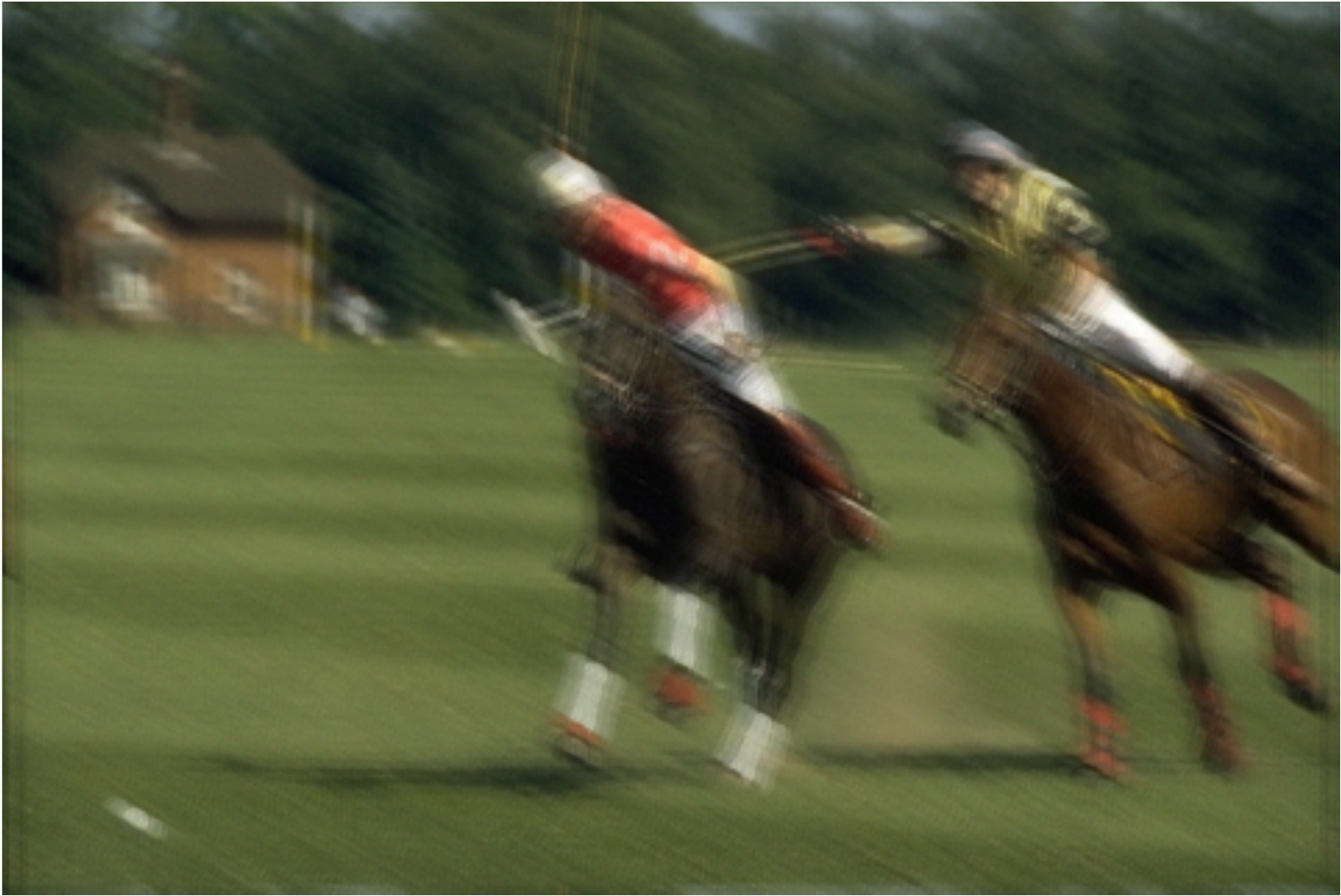}\\[.05em]
      \includegraphics[width=.495\linewidth,viewport=20 200 135 255,clip]{j_subimg1}\hfill%
      \includegraphics[width=.495\linewidth,viewport=270 160 385 215,clip]{j_subimg1}\\[-.7em]
    \end{minipage}%
  }\hspace{0.3cm}
  \subfigure[TViso {\scriptsize (27.10)}]{%
    \begin{minipage}{.470\linewidth}
    \includegraphics[width=1\linewidth]{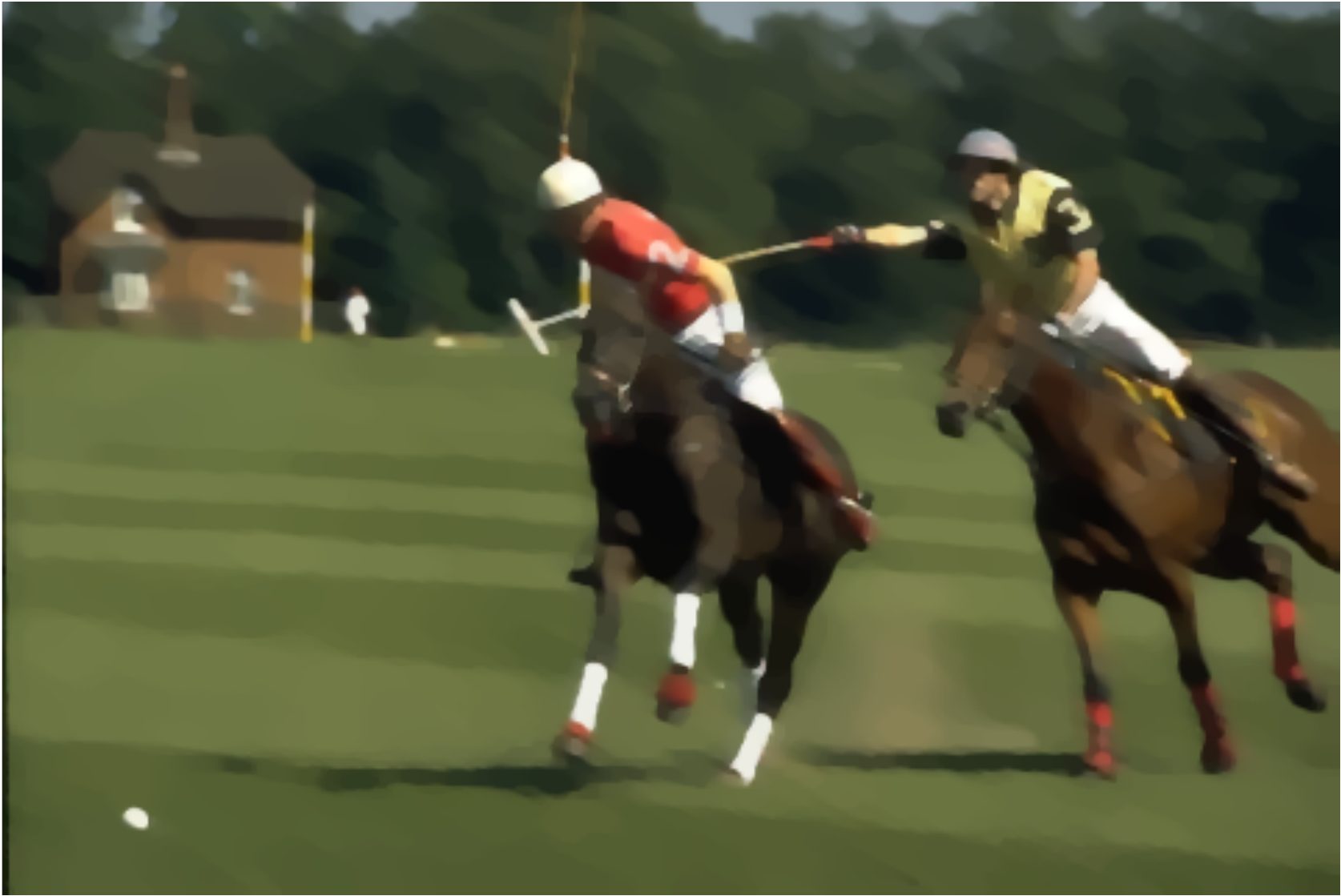}\\[.05em]
      \includegraphics[width=.495\linewidth,viewport=20 200 135 255,clip]{j_subimg2}\hfill%
      \includegraphics[width=.495\linewidth,viewport=270 160 385 215,clip]{j_subimg2}\\[-.7em]
    \end{minipage}%
  }\\[-.5em]
  \subfigure[QO {\scriptsize (29.57)}]{%
    \begin{minipage}{.470\linewidth}
      \includegraphics[width=1\linewidth]{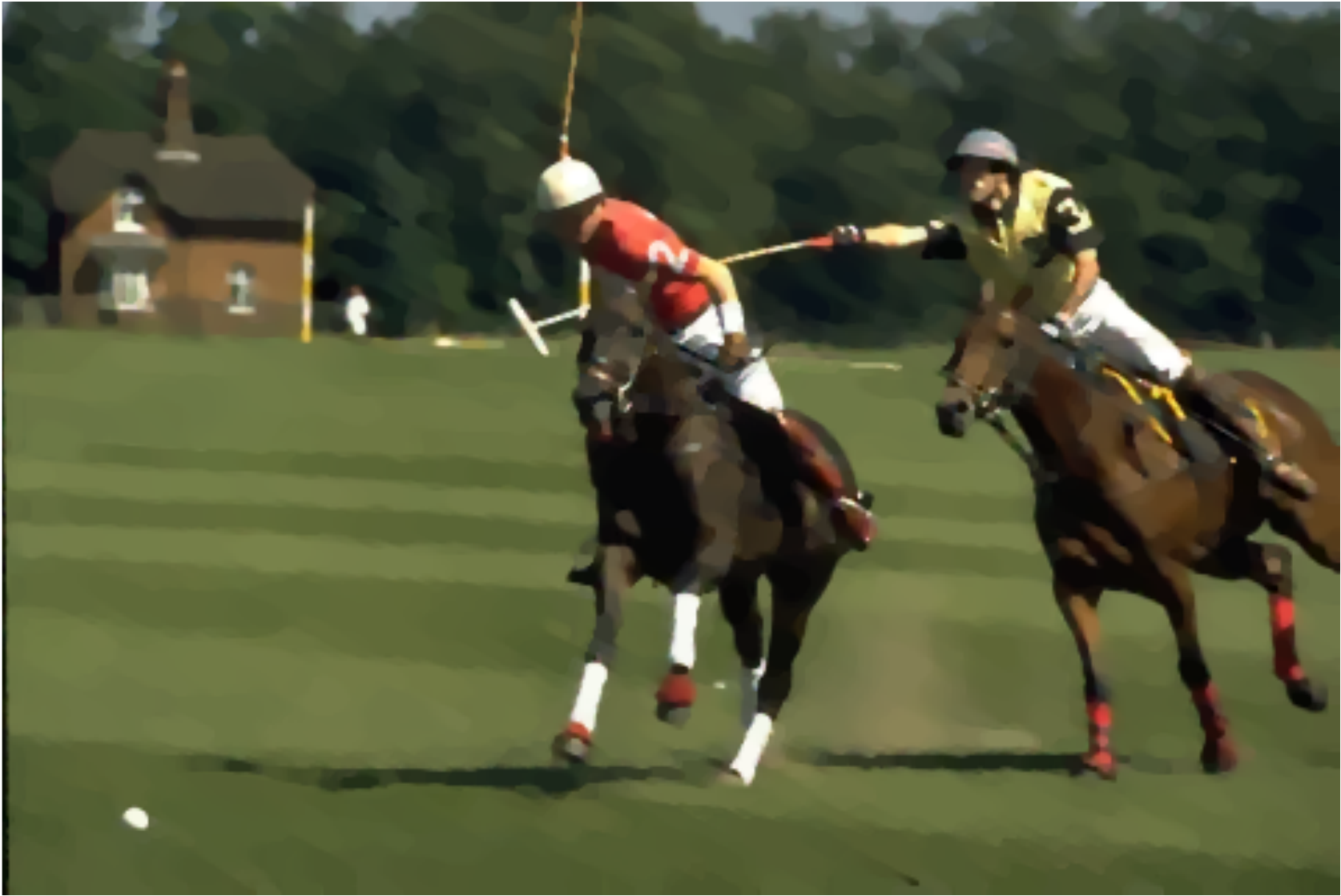}\\[.05em]
      \includegraphics[width=.495\linewidth,viewport=20 200 135 255,clip]{j_subimg5}\hfill%
      \includegraphics[width=.495\linewidth,viewport=270 160 385 215,clip]{j_subimg5}\\[-.7em]
    \end{minipage}%
  }\hspace{0.3cm}
  \subfigure[New SD {\scriptsize (30.35)}]{%
    \begin{minipage}{.470\linewidth}
      \includegraphics[width=1\linewidth]{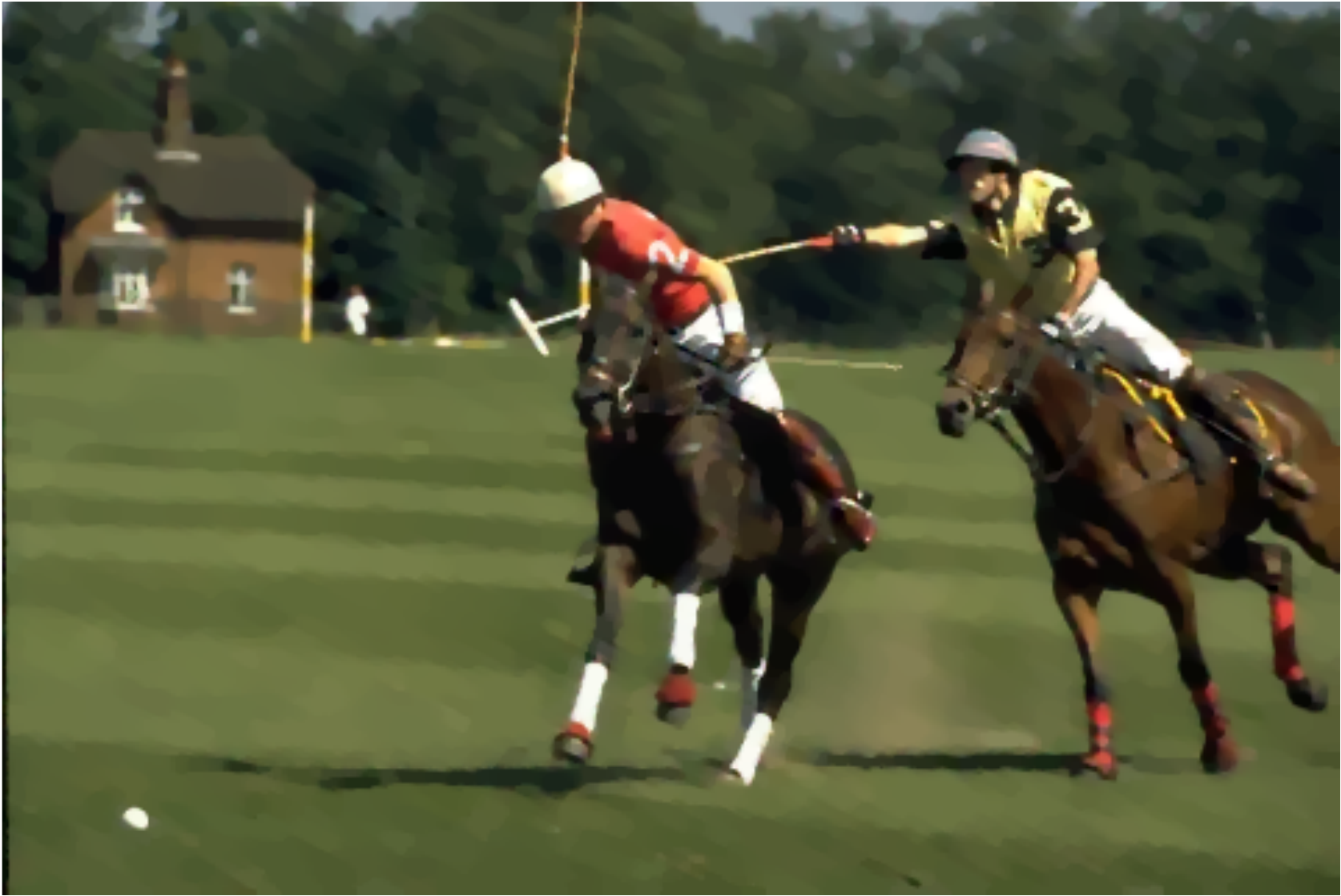}\\[.05em]
      \includegraphics[width=.495\linewidth,viewport=20 200 135 255,clip]{j_subimg6}\hfill%
      \includegraphics[width=.495\linewidth,viewport=270 160 385 215,clip]{j_subimg6}\\[-.7em]
    \end{minipage}%
  }\\
  \vspace{-1.2em}
  \caption{
    (a) An 8bit color image corrupted by a directional blur and
    Gaussian noise with standard deviation $\sigma=2$.
    (b) Solution of TViso.
    Debiased solution with (d) QO and (e) SD.
    The PSNR is indicated in brackets bellow each image.}
  \label{fig:deblurring}
  \vspace{-.em}
\end{figure}

We considered TViso regularization of degraded color images.
We defined blocks obtained by applying
$\Gamma=[\nabla_x^R,\nabla_y^R,\nabla_x^G,\nabla_y^G,\nabla_x^B,\nabla_y^B]$
where $m=n$, $b=6$, and $\nabla^C_d$ denotes the gradient in the direction $d \in \{ x, y \}$
for the color channel $C \in \{ R, G, B \}$.
We first focused on a denoising problem $y = x + w$ where
$x$ is an 8bit color image and $w$ is an additive white Gaussian
noise with standard deviation $\sigma=20$.
We next focused on a deblurring problem $y = \Phi x + w$ where
$x$ is an 8bit color image, $\Phi$ is a convolution
simulating a directional blur, and  $w$ is an additive white Gaussian
noise with standard deviation $\sigma=2$.
We chose $\lambda=4.3 \sigma$.
We applied the iterative primal-dual algorithm with our joint-refitting
(Algorithm \eqref{algo_final}) for $1,000$ iterations,
with $\tau = 1/4$, $\sigma = 1/6$ and $\theta = 1$.

Results are provided on Fig.~\ref{fig:denoising} and \ref{fig:deblurring}.
Comparisons of refitting with our proposed SD block penalty,
HD (only for denoising) and QO are provided.
Using our proposed SD block penalty offers the best refitting
performances in terms of both visual and quantitative measures.
The loss of contrast of TViso is well-corrected,
amplitudes are enhanced while smoothness and sharpness of TViso is
preserved. Meanwhile, the approach does not  create artifacts,
invert contrasts, or reintroduce information
that were not recovered by TViso.

\section{Conclusion}
In this work, we have reformulated the refitting problem of solutions promoted  by $\ell_{12}$ structured sparse analysis in terms of block penalties. We have introduced a new block penalty that interpolates between Bregman iterations \cite{Osher} and direction preservation \cite{brinkmann2016bias}.\footnotesize
This framework easily allows the inclusion of  additional desirable properties of refitted solutions as well as new  penalties that may increase refitting performances.
In order to take advantage of our efficient joint-refitting algorithm,  it is  important to consider  {\em simple} block penalty functions, which proximal operator   can be computed explicitly or at least easily.
Refitting in the case of other regularizers and  loss functions will be investigated in the future. \vspace{0.1cm}\\

\noindent
{
{\bf Acknowledgments}
This project has been carried out with  support from the French
State, managed by the French National Research Agency (ANR-16-CE33-0010-01).
This work was supported by the European Union’s Horizon 2020 research and innovation programme under the Marie Skłodowska-Curie grant agreement No 777826.\vspace{-0.1cm}
}

\appendix
\section{Proximity operators of block penalties}

\subsection{Convex conjugates $\phi^*$}
We  here compute the convex conjugate $\phi^*$ of the different  block penalties $\phi(z,\hat z)$ that only depends on $z\in \RR^b$ and where $\hat z\in \RR^b$ is a given fixed non null vector.
We consider the following representation of the vectors $z$ with respect to the $\hat z$ axis: $z=\alpha \frac{\hat z}{\norm{\hat z}}+\beta \frac{\hat z^\perp}{\norm{\hat z}}$.
This expression is valid for the case $b=2$. If $b=1$ then $z$ is only parameterized by $\alpha$. When $b>2$, $\hat z^\perp$ must be understood as a subspace $S$ of dimension $b-1$ and $\beta$ as a vector of $b-1$ components corresponding to each dimension of $S$.
With this change of variables, we have $\norm{z}^2=\alpha^2+\norm{\beta}^2$ (since $\beta$ is of dimension $b-1$),  $\cos(z,\hat z)=\alpha/\sqrt{\alpha^2+\norm{\beta}^2}$, $P_{\hat z}(z)=\alpha \hat z/\norm{\hat z}$ and $z-P_{\hat z}(z)=\beta \hat z^\perp/\norm{\hat z}$.
We also observe for instance that
   $|\cos(z,\hat z)|=1\Leftrightarrow \norm{\beta}=0$.
 All the block penalties $\phi(z, \hat z)$ can thus be expressed as   $\phi(\alpha,\beta)$.
The convex conjugate reads
\begin{equation}\label{convex_conj}
\phi^*(\alpha_0,\beta_0)=\sup_{\alpha,\beta} \alpha_0 \alpha+\dotp{\beta_0}{\beta}-\phi(\alpha,\beta).\\
\end{equation}
\paragraph{The block penalty $\phi_{\mathrm{HO}}$}
reads
$$\phi(z,\hat z)=\choice{
    0 & \ifq |\cos(z, \hat{z})| = 1~,\\
   + \infty & \otherwise~.
  }$$
  Hence, it gives
\begin{equation}\label{convex_conjHO}
\begin{split}
\phi^*(\alpha_0,\beta_0)&=\sup_{\alpha,\beta} \alpha_0 \alpha+\dotp{\beta_0}{\beta}-\choice{
    0 & \ifq \norm{\beta}=0~,\\
   + \infty & \otherwise}\\
    &=\choice{
    0 & \ifq \alpha_0=0 ~,\\
   + \infty & \otherwise~,
  }
\end{split}
\end{equation}
\paragraph{The block penalty $\phi_{\mathrm{HD}}$}
reads
$$\phi(z,\hat z)=\choice{
    0 & \ifq \cos(z, \hat{z}) = 1~,\\
  +  \infty & \otherwise~.
  }$$
  Hence, it gives
\begin{equation}\label{convex_conjHD}
\begin{split}
\phi^*(\alpha_0,\beta_0)&=\sup_{\alpha,\beta} \alpha_0 \alpha+\dotp{\beta_0}{\beta}-\choice{
    0 & \ifq \norm{\beta}=0\,\,\,\,\,\, \andq \alpha\geq 0~,\\
   + \infty & \otherwise}\\
    &=\choice{
    0 & \ifq \alpha_0\leq 0 ~,\\
   + \infty & \otherwise~,
  }
\end{split}
\end{equation}
\paragraph{The block penalty $\phi_{\mathrm{QO}}$}
reads
$$\phi(z,\hat z)=  \frac{\lambda}{2}
  \frac{\norm{z}^2}{\norm{\hat{z}}}(1-
  \cos^2(z, \hat{z})).$$
  Hence, it gives
\begin{equation}\label{convex_conjQO}
\begin{split}
\phi^*(\alpha_0,\beta_0)&=\sup_{\alpha,\beta} \alpha_0 \alpha+\dotp{\beta_0}{\beta}- \frac{\lambda}{2}
  \frac{\alpha^2+\norm{\beta}^2}{\norm{\hat{z}}}\left(1-
  \frac{\alpha^2}{\alpha^2+\norm{\beta}^2}\right)\\
  &=\sup_{\alpha,\beta} \alpha_0 \alpha+\dotp{\beta_0}{\beta}- \frac{\lambda}{2}
  \frac{\norm{\beta}^2}{\norm{\hat{z}}}\\
\end{split}
\end{equation}
The optimality condition on $\beta$ give
$$\beta=\frac{\beta_0\norm{\hat z}}{\lambda},$$
so that
$$\phi^*(\alpha_0,\beta_0)=\choice{
   \frac{\norm{\hat z}\norm{\beta_0}^2}{2\lambda} & \ifq \alpha_0\neq 0 ~,\\
   + \infty & \otherwise~.
}$$

\paragraph{The block penalty $\phi_{\mathrm{QD}}$}
reads
$$\phi(z,\hat z)=   \choice{ \frac{\lambda}{2}
  \frac{\norm{z}^2}{\norm{\hat{z}}}(1-
  \cos^2(z, \hat{z}))& \ifq \cos(z, \hat{z})\geq 0\\
     \frac{\lambda}{2}
  \frac{\norm{z}^2}{\norm{\hat{z}}}&\otherwise .
  }$$
  Hence, it gives
\begin{equation}\label{convex_conjQD}
\begin{split}
\phi^*(\alpha_0,\beta_0)&=\sup_{\alpha,\beta} \alpha_0 \alpha+\dotp{\beta_0}{\beta}-  \frac{\lambda}{2}
  \frac{\norm{\beta}^2}{\norm{\hat{z}}}-\choice{
 \frac{\lambda}{2} \frac{\alpha^2}{\norm{\hat{z}}} & \ifq \alpha\leq 0\\
0&\otherwise .
  }\\
\end{split}
\end{equation}
We observe that if $\alpha_0>0$, taking $\beta=0$ and letting $\alpha\to\infty$ leads to $\phi^*(\alpha_0,\beta_0)=+\infty$.
Next, the optimality conditions on $\beta$ and $\alpha$ give
$$\beta=\frac{\beta_0\norm{\hat z}}{\lambda},\hspace{1cm}\alpha=\frac{\alpha_0\norm{\hat z}}{\lambda},$$
so that
$$\phi^*(\alpha_0,\beta_0)=\choice{
   \frac{\norm{\hat z}(\alpha_0^2+\norm{\beta_0}^2)}{2\lambda} & \ifq \alpha_0\leq 0 ~,\\
   + \infty & \otherwise~.
}$$

\paragraph{The block penalty $\phi_{\mathrm{SD}}$}
reads
$\phi(z,\hat z)=  \frac\lambda2 \norm{z}(1-
  \cos(z, \hat{z})).
  $ 
  Hence, it gives
\begin{equation}\label{convex_conjSD}
\begin{split}
\phi^*(\alpha_0,\beta_0)&=\sup_{\alpha,\beta} \alpha_0 \alpha+\dotp{\beta_0}{\beta}- \frac\lambda2
  \sqrt{\alpha^2+\norm{\beta}^2 }\left(1-\frac{\alpha}{\sqrt{\alpha^2+\norm{\beta}^2}}\right)\\
  &=\sup_{\alpha,\beta} \alpha_0 \alpha+\dotp{\beta_0}{\beta}-  \frac\lambda2
\left(  \sqrt{\alpha^2+\norm{\beta}^2 }-\alpha\right)\\
  &=\sup_{\alpha,\beta} (\alpha_0+\lambda/2) \alpha+\dotp{\beta_0}{\beta}-  \frac\lambda2
\sqrt{\alpha^2+\norm{\beta}^2 }\\
\end{split}
\end{equation}
We observe that if $\sqrt{\norm{\beta_0}^2+(\alpha_0+\lambda/2)^2}>\lambda/2$, then  letting  $\alpha\to\sign{(\alpha_0+\lambda/2)}\times \infty$ and $\beta\to\sign{\beta_0}\times \infty$ leads to $\phi^*(\alpha_0,\beta_0)=+\infty$.
As a consequence we find
$$\phi^*(\alpha_0,\beta_0)=\choice{
  0& \ifq  \sqrt{\norm{\beta_0}^2+(\alpha_0+\lambda/2)^2}\leq\lambda/2~,\\
   + \infty & \otherwise~.
}$$

\subsection{Computing $\prox_{\sigma \phi^*}$}
We here give the computation of the proximal operator of the different $\phi^*$ that is given at point $(\alpha_0,\beta_0)$ by
\begin{equation}
  \prox_{\sigma\phi^*}(\alpha_0,\beta_0)=\uargmin{\alpha,\beta} \frac1{2\sigma}\left(\norm{\alpha-\alpha_0}^2+\norm{\beta-\beta_0}^2\right)+\phi^*(\alpha,\beta).
  \end{equation}
\paragraph{Block penalty $\phi_{\mathrm{HO}}$.} We have
\begin{equation}
\begin{split}\prox_{\sigma\phi^*}(\alpha_0,\beta_0)&=\uargmin{\alpha,\beta} \frac1{2\sigma}\left(\norm{\alpha-\alpha_0}^2+\norm{\beta-\beta_0}^2\right)+\choice{
    0 & \ifq \alpha=0 ~,\\
   + \infty & \otherwise~,
  }\\
  &=(0,\beta_0).
  \end{split}
\end{equation}
\paragraph{Block penalty $\phi_{\mathrm{HD}}$.} We have
\begin{equation}
\begin{split}\prox_{\sigma\phi^*}(\alpha_0,\beta_0)&=\uargmin{\alpha,\beta} \frac1{2\sigma}\left(\norm{\alpha-\alpha_0}^2+\norm{\beta-\beta_0}^2\right)+\choice{
    0 & \ifq \alpha\leq 0 ~,\\
   + \infty & \otherwise~,
  }\\
  &=(\min(0,\alpha_0),\beta_0).
  \end{split}
\end{equation}

\paragraph{Block penalty $\phi_{\mathrm{QO}}$.}
We have
\begin{equation}
\begin{split}\prox_{\sigma\phi^*}(\alpha_0,\beta_0)&=\uargmin{\alpha,\beta} \frac1{2\sigma}\left(\norm{\alpha-\alpha_0}^2+\norm{\beta-\beta_0}^2\right)+\choice{
   \frac{\norm{\hat z}\norm{\beta}^2}{2\lambda} & \ifq \alpha\neq 0 ~,\\
   + \infty & \otherwise~.
}\\
  &=\frac{\lambda}{\lambda+\sigma \norm{\hat z}}\left(0,\beta_0\right),
  \end{split}
\end{equation}
since the optimality condition with respect to $\beta$ gives
$\lambda(\beta-\beta_0)+\sigma \norm{\hat z} \beta=0$.

\paragraph{Block penalty $\phi_{\mathrm{QD}}$.} We have
\begin{equation}
\begin{split}\prox_{\sigma\phi^*}(\alpha_0,\beta_0)&=\uargmin{\alpha,\beta} \frac1{2\sigma}\left(\norm{\alpha-\alpha_0}^2+\norm{\beta-\beta_0}^2\right)+\choice{
   \frac{\norm{\hat z}(\alpha^2+\norm{\beta}^2)}{2\lambda} & \ifq \alpha\leq 0 ~,\\
   + \infty & \otherwise~.
}\\
  &=\frac{\lambda}{\lambda+\sigma \norm{\hat z}}\left(\min(0,\alpha_0),\beta_0\right).
  \end{split}
\end{equation}

\paragraph{Block penalty $\phi_{\mathrm{SD}}$.} We have
\begin{equation}
\begin{split}&\prox_{\sigma\phi^*}(\alpha_0,\beta_0)\\=&\uargmin{\alpha,\beta} \frac1{2\sigma}\left(\norm{\alpha-\alpha_0}^2+\norm{\beta-\beta_0}^2\right)+\choice{
  0& \ifq  \sqrt{\norm{\beta}^2+(\alpha+\lambda/2)^2}\leq\lambda/2~,\\
   + \infty & \otherwise~.
}\\
  = &  \frac\lambda{2}\frac{(\alpha_0+\lambda/2,\beta_0) }{\max(\lambda/2,\sqrt{\norm{\beta_0}^2+(\alpha_0+\lambda/2)^2})} -\left(\lambda/{2},0\right),  
  \end{split}
\end{equation}
which just corresponds to the projection of the $\ell_2$ ball of $\RR^b$ of radius $\lambda/2$ and center $(-\lambda/2,0)$.

\bibliographystyle{abbrv}

\end{document}